\documentclass[11pt]{article}
\usepackage[T1]{fontenc} 
\usepackage[utf8]{inputenc}

\usepackage{lmodern}

\usepackage{secdot}
\usepackage{amssymb}
\usepackage{amsmath}
\usepackage{latexsym}
\usepackage{calrsfs}
\usepackage{graphicx,color,psfrag}
\usepackage{natbib}

\newtheorem{theorem}{Theorem}[section]
\newtheorem{lemma}[theorem]{Lemma}

\newtheorem{proposition}[theorem]{Proposition}
\newtheorem{definition}[theorem]{Definition}
\newtheorem{assumption}[theorem]{Assumption}
\newenvironment{proof}{\trivlist\item[]\emph{Proof}.}
{\unskip\nobreak\hskip 1em plus 1fil\nobreak$\Box$
\parfillskip=0pt\endtrivlist}

\def\argmax{\mathop{\rm arg\,max}}

\newcommand{\Ex}{\ensuremath{\mathsf{E}}}

\usepackage[colorlinks=true,urlcolor=blue,citecolor=blue,linkcolor=blue]{hyperref}
\usepackage[noabbrev,nameinlink]{cleveref}     

\def\EMAIL#1{\href{mailto:#1}{#1}}

\begin{document}

\title{A faster index algorithm and a computational study  for bandits with switching costs}


\author{Jos\'e Ni\~no-Mora 
\\ Department of Statistics \\
    Carlos III University of Madrid \\
     28903 Getafe (Madrid), Spain \\  \EMAIL{jose.nino@uc3m.es}, \href{http://alum.mit.edu/www/jnimora}{http://alum.mit.edu/www/jnimora} \\
      ORCID: \href{http://orcid.org/0000-0002-2172-3983}{0000-0002-2172-3983}}
 
\date{Published in \textit{INFORMS Journal on Computing}, vol.\  20, pp. 255–269,  2008 \\ \vspace{.1in}
DOI: \href{https://doi.org/10.1287/ijoc.1070.0238}{10.1287/ijoc.1070.0238}}

\maketitle


\begin{abstract}%
We address the intractable 
 multi-armed bandit problem
with  switching costs, for which
Asawa and Teneketzis introduced in  [M.\ Asawa and D.\ Teneketzis. 1996. Multi-armed bandits with switching
 penalties. IEEE Trans. Automat. Control, 41 328--348]  an
 index that
partially characterizes optimal policies, 
attaching
 to each project state
 a ``continuation index'' (its Gittins index) and a ``switching index.''
They proposed to jointly compute both as the
Gittins index of a project with $2 n$
states --- when the original project has $n$ states --- resulting in an eight-fold increase in $O(n^3)$ arithmetic
operations
relative to those to compute the continuation index.
We present a faster decoupled 
computation method, which in a first stage computes the 
continuation index 
 and then, in a second stage, computes the switching index
an order of magnitude faster in at most $n^2 + O(n)$ arithmetic operations,
 achieving overall a four-fold reduction in arithmetic operations
and substantially reduced memory operations.
The analysis exploits the fact that the Asawa and Teneketzis
index is the marginal productivity index of
the project in its
restless reformulation, using methods 
 introduced by the author. 
Extensive computational experiments are reported, which
demonstrate the dramatic runtime speedups achieved by the new
 algorithm, as well as
the near-optimality of the resultant index policy and its substantial
gains
against the benchmark Gittins index policy across a wide 
range of randomly generated two- and three-project instances. 
\end{abstract}%


\newpage

\section{Introduction}
\label{s:intro}
The \emph{Gittins index} furnishes a
tractable solution to the
\emph{multi-armed bandit problem} (MABP).
This involves finding an optimal sequential effort allocation 
policy for a finite collection of stochastic projects, one of which
must be engaged at each discrete time period over an infinite
horizon. 
Projects are modeled as
\emph{bandits}, i.e., binary-action (1: active; 0: passive) 
 \emph{Markov decision processes} (MDPs) than can only change state
when active.
The optimal policy engages at each
time a project of largest index, where the index
 is defined separately for each project as a
function of its state.
See 
 \cite{gijo74} and  \cite{gi79}.

A critical assumption underlying such a result
 is that switching projects is costless.
Yet, as noted in \citet{banksun94}, ``it is difficult to imagine
a relevant economic decision problem in which the decision-maker may
costlessly move between alternatives.''
Incorporation of such costs yields the 
\emph{multi-armed bandit problem with 
switching costs} (MABPSC). See, e.g.,
\citet{aggraetal88}, \citet{votenek94}, and the survey of
 \citet{jun04}.

In light of the Gittins and Jones result for the MABP, 
it is appealing to try to design good index policies 
for the MABPSC. 
Any such a policy  attaches 
 an index $\nu_m(a_m^-, i_m)$ to each project
$m$, which is a function of its previous
action $a_m^-$ and current state $i_m$, thus
decoupling into a 
``continuation index'' $\nu_m(1, i_m)$ and a 
``switching index'' $\nu_m(0, i_m)$.
The index policy prescribes to engage at each time
a project of currently largest index.

Though index policies are generally 
suboptimal for the MABPSC (cf.\ 
 \cite{banksun94}),  \cite{asatene96} introduced
 an index that
 partially characterizes optimal policies.
Their continuation index is the 
 Gittins index, while
their switching index is the maximum rate of 
expected discounted reward  minus switching cost 
 per unit of 
expected discounted time achievable by stopping rules that
engage an initially passive project.

\cite{asatene96} proposed to jointly compute both indices by: (i) formulating a modified
 project 
\emph{without} switching costs, yet having \emph{twice}
 the number of states --- the  $(a_m^-, i_m)$;
 and (ii) computing the Gittins index of the latter. 
While several algorithms are available to compute the Gittins
index, at present the best complexity for a general
$n$-state project is achieved by the 
 \emph{fast-pivoting} algorithm in 
\cite{nmijoc06}, which performs $(2/3) n^3 + O(n^2)$ 
arithmetic operations. 
Yet using it to compute the
Asawa and Teneketzis (AT) indices via their scheme yields
an operation count of $(2/3) (2n)^3 + O(n^2)$: an  eight-fold
increase relative to the effort to compute the continuation index
 alone, 
which severely hinders deployment  to
large-scale models.
This raises the need to develop a significantly faster
computation method, which is the prime goal
of this paper.

We accomplish such a goal via a seemingly indirect
route, which exploits the reformulation of a classic bandit
with switching costs as a \emph{restless bandit} 
---  one that can 
change state when passive --- \emph{without}
switching costs.
This allows us to deploy the approach to restless bandit 
indexation 
introduced by 
 \cite{whit88} and developed by the author (cf.\ \citet{nmtop07}),
which defines a \emph{marginal productivity index} (MPI) for a restricted
class of bandits termed \emph{indexable}.
This forms the basis for 
a heuristic index policy for the intractable
\emph{multi-armed restless bandit problem} (MARBP). 

We prove and exploit the fact
 that, for bandits with switching costs,
 their MPI and their AT
index are  the same.
This provides
an economically insightful justification for use of such an
index policy as a heuristic for the MABPSC, based
on the economic characterization given in \citet{nmmp02,nmmor06} of the 
MPI as a state-dependent measure of the marginal
 productivity of work on a project.
MPI-based priority policies dynamically allocate effort 
where it appears to be currently more
productive, using the MPI of each project as a proxy --- as it 
ignores interactions --- marginal productivity
measure. 

While  \cite{asatene96} focused
on projects with
constant startup costs, we allow
state-dependent startup and shutdown costs, provided
 that their sum  is
nonnegative at each state.
This ensures that the continuation index is 
at least as large as the switching index, consistently 
with the intuitive \emph{hysteresis} property of optimal policies
noticed in \cite{banksun94}: 
``it is obvious that in comparing two otherwise
identical arms, one of which was used in the previous period, the
one which was in use must necessarily be more attractive than the 
one which was idle.''. See also
\citet{duhon03}.

We deploy the methods introduced in 
\citet{nmaap01,nmmp02,nmmor06}, showing that the restless bandits of concern are
\emph{PCL-indexable} --- after their
satisfaction of \emph{partial conservation laws} (PCLs), which allows us to compute the  MPI using the
 \emph{adaptive-greedy algorithm} introduced in such work.
We further decouple the latter into
a faster two-stage method: a first stage that computes the continuation
 index along with additional quantities; and a
 second stage 
that uses the first-stage's output to compute the switching index.

To implement such a scheme, one can use for the first stage several
algorithms in \cite{nmijoc06}, such as the \emph{fast-pivoting algorithm with extended output}
 FP(1), which 
performs $(4/3) n^3 + O(n^2)$ arithmetic operations.
For the second stage,  we present a 
switching-index algorithm that performs \emph{at most} $n^2 + O(n)$ arithmetic
operations, rendering negligible
the marginal effort to compute the switching index.
Relative to the
 AT scheme, our two-stage method
achieves overall a four-fold reduction 
in arithmetic operations along with substantially reduced memory
operations.
Such an algorithm is the main contribution of this paper.

A computational study
demonstrates that such an improved complexity 
yields dramatic runtime speedups.
The study is complemented by a set of experiments that 
demonstrate the near-optimality of the index policy and its substantial
gains
against the benchmark Gittins-index policy across an extensive
range of two- and three-project instances.

Section \ref{s:MBPsc} discusses the model and its restless-bandit
reformulation.
Section \ref{s:impipcl} reviews restless bandit indexation, gives an indexability proof based on
dynamic programming (DP), shows that the MPI is precisely the
AT index, and outlines
the approach to be deployed.
Section \ref{s:pclapsc} carries out a PCL-indexability 
analysis of a project with switching costs, in its restless reformulation. 
Section \ref{s:ecsi} draws on such an analysis to develop the new
decoupled MPI algorithm. 
Section \ref{s:es} reports the computational
 study's results.
Section \ref{s:concl} concludes.

\section{MABPSC and Restless-Bandit Reformulation}
\label{s:MBPsc}
\subsection{Model Formulation}
\label{s:tm}
We next describe the MABPSC outlined in Section \ref{s:intro}.
Consider a collection of $M$ finite-state  projects,
one of which must be engaged (\emph{active})
at each  discrete time period $t =0, 1, 2, \ldots$ over an infinite
horizon, while the others are rested (\emph{passive}).
When project $m$ occupies state $i_m$ --- belonging in its state space
$N_m$ --- and is engaged, it yields an
 \emph{active reward}
$R_m^1(i_m) = R_m(i_m)$ and its state moves
to $j_m$  with  probability $p_m(i_m, j_m)$.
If the project is instead rested, it yields a zero
\emph{passive reward} $R_m^0(i_m) = 0$ and its state remains frozen.

Switching projects is costly.
When project $m$ occupies state $i_m$ and is freshly engaged
(resp.\ rested), a 
\emph{startup cost} $c_m(i_m)$
(resp.\ \emph{shutdown cost} $d_m(i_m)$) is incurred. 
We assume that
$c_m(i_m) + d_m(i_m) \geq 0$.
Rewards and costs are time-discounted with factor $0 < \beta
< 1$.

Actions are chosen by adoption of a \emph{scheduling policy}
$\boldsymbol{\pi}$, drawn from the class $\boldsymbol{\Pi}$ of \emph{admissible policies}, which
are nonanticipative relative to the  history 
of states and actions and engage one project at a time.
The  MABPSC
is to find an admissible
policy
maximizing the expected total
discounted value of rewards earned minus switching costs incurred.

Denote by $X_m(t) \in N_m$ and
$a_m(t) \in \{0, 1\}$ the 
state and action for project $m$ at time $t$, respectively.
We will use the notation 
\begin{equation}
\label{eq:amminus}
a_m^-(t) \triangleq a_m(t-1), \quad
\bar{a}_m(t) \triangleq 1 - a_m(t), \quad
\text{and} \quad \bar{a}_m^-(t) \triangleq 
 \bar{a}_m(t-1).
\end{equation}
Since initial conditions must specify 
whether each project $m$ is initially set up, 
we denote such a status by $a_m^-(0)$.
We further define the project's 
\emph{augmented state} $\hat{X}_m(t) \triangleq (a_m^-(t),
X_m(t))$, which moves over the \emph{augmented state space}
$\hat{N}_m \triangleq \{0, 1\} \times N_m$. 
The system's \emph{joint augmented state} is thus $\hat{\mathbf{X}}(t)
\triangleq (\hat{X}_m(t))_{m=1}^M$ and its \emph{joint action process}
is $\mathbf{a}(t) \triangleq (a_m(t))_{m=1}^M$.

Given a joint initial  state 
$\hat{\mathbf{X}}(0) =
\hat{\boldsymbol{\imath}} = ((a_m^-, i_m))_{m=1}^M$
 we  can  formulate the MABPSC as
\begin{equation}
\label{eq:vui}
\begin{split}
& \max_{\boldsymbol{\pi} \in \boldsymbol{\Pi}}
  \Ex_{\boldsymbol{\hat{\imath}}}^{\boldsymbol{\pi}}\left[\sum_{m=1}^M \sum_{t=0}^\infty   
  \big\{R_m^{a_m(t)}\left({X_m(t)}\right) - 
  c_m\left({X_m(t)}\right) \bar{a}_m^-(t) a_m(t)  - 
  d_m\left({X_m(t)}\right) a_m^-(t) \bar{a}_m(t)\big\}
  \beta^t\right],
\end{split}
\end{equation}
where $\Ex_{\boldsymbol{\hat{\imath}}}^{\boldsymbol{\pi}}\left[\cdot\right]$ is
 expectation under $\boldsymbol{\pi}$ conditioned on the initial joint state
 being equal to $\boldsymbol{\hat{\imath}}$.

\subsection{Restless-Bandit Reformulation}
\label{s:rbref}
Taking $\hat{X}_m(t)$ as the state of project $m$ 
 yields a reformulation of  (\ref{eq:vui}) 
as an MARBP (cf.\ Section \ref{s:intro})
\emph{without} switching costs.
The rewards and transition probabilities in such a  reformulation are as follows.
If the project occupies augmented state $(a_m^-,
i_m)$ and is engaged, the 
active reward $\hat{R}_m^1(a_m^-, i_m) \triangleq 
R_m^1(i_m) - c_m(i_m) \bar{a}_m^-$ accrues, and its augmented state
moves to 
$(1, j_m)$ with active transition probability
$\hat{p}_m^1\big((a_m^-, i_m), (1, j_m)\big) \triangleq p_m(i_m,
j_m)$.
If it is rested, the one-period
passive reward $\hat{R}_m^0(a_m^-, i_m) \triangleq 
R_m^0(i_m) - d_m(i_m) a_m^-$ accrues, and its augmented state moves to
$(0, i_m)$ with a unity passive transition probability:
$\hat{p}_m^0\big((a_m^-, i_m), (0, i_m)\big) \equiv 1$.

We can thus reformulate (\ref{eq:vui}) as the 
MARBP
\begin{equation}
\label{eq:valref}
\max_{\boldsymbol{\pi} \in \boldsymbol{\Pi}}
\Ex_{\hat{\boldsymbol{\imath}}}^{\boldsymbol{\pi}}\left[\sum_{m=1}^M  \sum_{t=0}^\infty
  \hat{R}_{m}^{a_m(t)}\big({\hat{X}_m(t)}\big) \beta^t\right].
\end{equation}
This allows us to deploy 
the approach to 
restless bandit indexation introduced
by \cite{whit88} and developed by the author in recent work (cf.\ \citet{nmtop07}).
It must be
  emphasized, though, that the resultant index policies are generally suboptimal.

\section{Restless Bandit Indexation}
\label{s:impipcl}
In this section we review restless bandit indexation, give an indexability proof based on
DP, show that the MPI is the
AT index, and outline
the approach to be deployed.
We draw on  
\cite{whit88} and  \cite{nmaap01,nmmp02,nmmor06} for
the definition, analysis and
computation of the MPI, as it applies to a single project
as above --- in its restless reformulation.
We thus drop below the project label $m$ so that now, 
 e.g.,  $\hat{N} \triangleq \{0, 1\}
\times N$ denotes the project's augmented state space. 
We denote by $\pi$ and $\Pi$ the policies and policy space for
operating the project,  where the notation distinguishes them from
their boldface counterparts $\boldsymbol{\pi}$ and $\boldsymbol{\Pi}$
used above in the multi-project setting. 

We assume henceforth
that switching costs satisfy the following key condition.
\begin{assumption}
\label{ass:key}
$c_i + d_i \geq 0$ for $i \in N$.
\end{assumption}

\subsection{Indexability,  MPI, and Hysteresis}
\label{s:impi}
We
 evaluate the value of \emph{net rewards} (i.e., accounting for switching
 costs) earned on the project
under a policy $\pi \in
\Pi$ starting at  $(a_0^-, i_0) \in \hat{N}$ by the
\emph{reward measure}
\[
f^{\pi}_{(a_0^-, i_0)} \triangleq \Ex_{(a_0^-, i_0)}^{\pi}\left[\sum_{t=0}^\infty 
   \hat{R}^{a(t)}\big(\hat{X}(t)\big) \beta^t\right],
\]
and further evaluate the corresponding amount of work expended by the
\emph{work measure}
\[
g^{\pi}_{(a_0^-, i_0)} \triangleq \Ex_{(a_0^-, i_0)}^{\pi}\left[\sum_{t=0}^\infty 
  a(t) \beta^t\right].
\]

Imagining that work on the project
is paid for at the  \emph{wage} rate 
$\nu$ per unit work performed
leads us to consider the \emph{$\nu$-wage problem}
\begin{equation}
\label{eq:spnuwp}
\max_{\pi \in \Pi} f^{\pi}_{(a_0^-, i_0)} -
 \nu g^{\pi}_{(a_0^-, i_0)},
\end{equation}
which is to find an admissible policy achieving the maximum value
$v^{*}_{(a_0^-, i_0)}(\nu)$ of 
net rewards earned minus labor costs incurred.
We will use (\ref{eq:spnuwp}) to \emph{calibrate} 
the \emph{marginal value of work} at each state, by analyzing the
structure of its optimal
policies as the wage $\nu$ varies.

Since (\ref{eq:spnuwp}) is a finite-state and -action
MDP,  we know that for any wage $\nu \in
\mathbb{R}$ there exists an optimal policy that is 
stationary deterministic and
independent of the initial  state. Each such policy is
characterized by its \emph{active-state set}, or \emph{active set}, which is the subset of augmented states where
it prescribes to engage the project.
We write active sets as
\[
S_0 \oplus S_1
  \triangleq \{0\} \times S_0 \cup \{1\}
  \times S_1, \quad S_0, S_1 \subseteq N.
\] 
Thus,   the policy that
 we denote by $S_0 \oplus S_1$
engages the project when it was previously rested (resp.\ engaged) 
if its original 
state $X(t)$ lies in $S_0$ (resp.\ in $S_1$).

Hence, for any $\nu$ there is a unique 
\emph{minimal optimal active set} $S_0^*(\nu) \oplus S_1^*(\nu)
\subseteq \hat{N}$, consisting of all augmented states where
engaging the project is the only optimal action.
The dependence of such sets on
  $\nu$ is used next to define the \emph{indexability}
property.

\begin{definition}[Indexability; MPI]
\label{def:indx}
\textup{We say that the project is \emph{indexable}  if
  there exists an \emph{index} $\nu^*_{(a^-, i)}$ for 
 $(a^-, i) \in \hat{N}$ such that 
\begin{equation}
\label{eq:dstarnu}
S_0^*(\nu) 
 = \big\{(0, i)\colon \nu^*_{(0, i)} > \nu\big\}
\quad \text{and} \quad
S_1^*(\nu) 
 = \big\{(1, i)\colon \nu^*_{(1, i)} > \nu\big\}, \quad
\nu \in \mathbb{R}.
\end{equation}
We say that $\nu^*_{(a^-, i)}$ is the project's 
\emph{marginal productivity index}  (MPI), terming
$\nu^*_{(1, i)}$ the \emph{continuation MPI} and 
$\nu^*_{(0, i)}$ the \emph{switching MPI}.
}
\end{definition}

Thus, the project  is indexable with MPI 
$\nu^*_{(a^-, i)}$ if
it is optimal
to engage (resp.\ rest) the project when it occupies 
$(a^-, i)$ iff $\nu^*_{(a^-, i)} \geq \nu$ (resp.\ 
$\nu^*_{(a^-, i)} \leq \nu$).
\citet{whit88} introduced the concept of indexability for restless
bandits under the average criterion.
We introduced the term MPI in \citet{nmmor06},
motivated by the characterization of the index   as a
measure of the state-dependent
marginal productivity of work on a project.

To establish indexability and compute the MPI, 
\citet{nmaap01,nmmp02,nmmor06} develops an approach based on
guessing  the structure of optimal active
sets, as an \emph{active-set family} 
$\hat{\mathcal{F}} \subseteq 2^{\hat{N}}$ that contains sets
$S_0^*(\nu) \oplus S_1^*(\nu)$ as $\nu$ varies.
The intuition that, under Assumption \ref{ass:key}, optimal policies
should have the hysteretic
property (cf.\ Section \ref{s:intro})
that, if it is optimal to engage a project when it was previously
rested then
it should also be optimal to engage it when it was previously
active, leads us to posit that the right choice of $\hat{\mathcal{F}}$
must be
\begin{equation}
\label{eq:fhat}
\hat{\mathcal{F}} \triangleq \big\{S_0 \oplus S_1\colon S_0
  \subseteq S_1 \subseteq N\big\}.
\end{equation}

\begin{definition}[$\hat{\mathcal{F}}$-indexability]
\label{def:findx}
\textup{We say that the project is \emph{$\hat{\mathcal{F}}$-indexable} if:
(i) it is indexable; and
(ii) $S_0^*(\nu) \oplus S_1^*(\nu)
 \in \hat{\mathcal{F}}$ for $\nu \in \mathbb{R}$.
}
\end{definition}

Note that $\hat{\mathcal{F}}$ represents a family of policies. We will establish in Theorem \ref{the:indxb} that, 
with $\hat{\mathcal{F}}$ defined as in
(\ref{eq:fhat}), the restless projects of concern in this paper
are indeed $\hat{\mathcal{F}}$-indexable.

\subsection{Reduction to the Normalized Zero Shutdown Costs Case}
\label{s:rzscc}
This section shows that it suffices to consider
the zero shutdown costs case, extending a result in
\citet[Sec.\ 3]{banksun94} for the constant switching costs case $c_i
\equiv c, d_i \equiv d$.
Suppose that at a certain time, which we take to be $t = 0$,
 a project is freshly engaged for a random
duration given by a stopping-time rule $\tau$.
Denoting by $\mathbf{R} = (R_j)$,
$\mathbf{c} = (c_j)$ and $\mathbf{d} = (d_j)$ its
state-dependent active
reward, startup and shutdown cost vectors, we can write the
expected discounted net earnings over such a time span starting at
$X(0) = i$ as
\begin{equation}
\label{eq:frcd}
\begin{split}
f_i^\tau(\mathbf{R}, \mathbf{c}, \mathbf{d}) & \triangleq 
 \Ex_i^\tau\left[-c_i + \sum_{t = 0}^{\tau-1} R_{X(t)} \beta^t -
 d_{X(\tau)} \beta^\tau\right].
\end{split}
\end{equation}

We have the following result, where $\mathbf{I}$ is the identity
matrix indexed by the state space $N$, $\mathbf{P} = (p_{ij})_{i, j
  \in N}$ is the transition probability matrix, and $\mathbf{0}$ is a
vector of zeros.
\begin{lemma}
\label{lma:osgi}
$f_i^\tau(\mathbf{R}, \mathbf{c}, \mathbf{d}) = 
f_i^\tau\big(\mathbf{R} + (\mathbf{I} - \beta \mathbf{P}) \mathbf{d},
\mathbf{c} + \mathbf{d}, \mathbf{0}\big)$.
\end{lemma}
\begin{proof}
Use the elementary identity 
\[
d_{X(\tau)} \beta^\tau = d_i - \sum_{t=0}^{\tau-1} \{d_{X(t)} - \beta 
 d_{X(t+1)}\} \beta^t
\]
to obtain
\begin{align*}
f_i^\tau(\mathbf{R}, \mathbf{c}, \mathbf{d}) & 
\triangleq -c_i + \Ex_i^\tau\left[\sum_{t=0}^{\tau-1} R_{X(t)} \beta^t - d_{X(\tau)} 
  \beta^{\tau}\right] \\
& = -c_i -d_i + \Ex_i^\tau\left[\sum_{t=0}^{\tau-1}
    \big\{R_{X(t)} + d_{X(t)} - \beta d_{X(t+1)}\big\}
    \beta^t\right] = f_i^\tau\big(\mathbf{R} + (\mathbf{I} - \beta \mathbf{P}) \mathbf{d},
\mathbf{c} + \mathbf{d}, \mathbf{0}\big).
\end{align*}
\end{proof}

Lemma \ref{lma:osgi} shows how to eliminate shutdown costs: one need simply incorporate them into modified
startup costs and
active rewards given by the transformations
\begin{equation}
\label{eq:chatrhat}
\tilde{\mathbf{c}} \triangleq \mathbf{c} +
\mathbf{d} \quad \text{and} \quad 
\tilde{\mathbf{R}} \triangleq \mathbf{R} + (\mathbf{I} -
\beta \mathbf{P}) \mathbf{d}.
\end{equation} 
Note that in the case $c_j \equiv c$ and
$d_j \equiv d$ discussed in \citet{banksun94} such transformations
reduce to
$\tilde{c}_j \equiv c + d$ and 
$\tilde{R}_j = R_j + (1-\beta) d$, in agreement with their results.

For simplicity, we will hence focus the discussion henceforth in
the \emph{normalized} case $c_i \geq 0$, $d_i \equiv
0$, assuming that the transformations in (\ref{eq:chatrhat}) have been
carried out if required.

\subsection{AT Index, MPI, and 
 Indexability Proof}
\label{s:dpaidmpi}
We proceed to give representations of the
continuation and switching AT indices in the normalized case, 
 to present a DP-based proof of
$\hat{\mathcal{F}}$-indexability, and to show that the AT indices
match their MPI counterparts. 
The representations are given in terms of 
work and reward measures $g_i^S$ and $f_i^S$ for the underlying
 project without switching costs starting at $i \in N$, under the policy that
 engages the project when its original state $X(t)$ lies in $S \subseteq N$. 

Define the continuation and the switching AT index  by 
\begin{equation}
\label{eq:nu1i}
\nu_{(a^-, i)}^{\textup{AT}} \triangleq \max_{\tau > 0}
\frac{\displaystyle \Ex_i^\tau\left[-(1-a^-) c_i + \sum_{t=0}^{\tau-1} R_{X(t)}
    \beta^t\right]}{\displaystyle \Ex_i^\tau\left[\sum_{t=0}^{\tau-1}
    \beta^t\right]} = \max_{S \subseteq N\colon i \in S}
\frac{f^{S}_{i} - (1-a^-) c_i}{g^{S}_{i}},
\end{equation}
 where
$\tau$ denotes a stopping-time rule. 
Hence $\nu_{(1, i)}^{\textup{AT}}$ is precisely the project's Gittins index. 
Note that the right-most identity in
 (\ref{eq:nu1i}) follows since it suffices to restrict
 attention to stationary deterministic stopping rules, which are represented by
their continuation sets.

The next result shows that such indices are consistent with the 
policy family $\hat{\mathcal{F}}$ in (\ref{eq:fhat}).
\begin{lemma}
\label{lma:nu1gtnu0}
$\nu_{(1, i)}^{\textup{AT}} \geq \nu_{(0, i)}^{\textup{AT}}$, for $i \in N$.
\end{lemma}
\begin{proof}
Fix $i \in N$. For each set $S \subseteq N$,
Assumption \ref{ass:key}, i.e., $c_i \geq 0$, implies
\[
\frac{f_i^{S}}{g_i^{S}} \geq 
\frac{f^{S}_{i} - c_i}{g^{S}_{i}}.
\]
The result follows by maximizing over $S \subseteq N$ in each side of
latter inequality, which preserves it, and using the definitions of 
$\nu_{(1, i)}^{\textup{AT}}$ and $\nu_{(0, i)}^{\textup{AT}}$ in (\ref{eq:nu1i}).
\end{proof}

We next set out to use DP to show that the indices defined in
(\ref{eq:nu1i}) are indeed the project's continuation and
switching MPIs in Definition \ref{def:indx}.
We will use the \emph{Bellman equations} on the optimal value function $v_{(a^-,
  i)}^*(\nu)$ for $\nu$-wage problem
(\ref{eq:spnuwp}): for $i \in N$,
\begin{equation}
\label{eq:be}
v_{(a^-, i)}^*(\nu) = \max \big\{\beta v_{(0, i)}^*(\nu), R_i  - (1-
a^-) c_i
 - \nu + \beta 
 \sum_{j \in N} p_{ij} v_{(1, j)}^*(\nu)\big\}.
\end{equation}

We start by establishing that it suffices to restrict
attention to active-set family $\hat{\mathcal{F}}$.

\begin{lemma}
\label{lma:intuit}
If it is optimal to rest the project in $(1, i)$, then it is
optimal to rest it in  $(0, i)$.
\end{lemma}
\begin{proof}
The result follows from the implication
\[
\beta v_{(0, i)}^*(\nu) \geq R_i - \nu + \beta \sum_{j \in N}
 p_{ij} v_{(1, j)}^*(\nu) \Longrightarrow \beta v_{(0, i)}^*(\nu) \geq -c_i 
  +  R_i - \nu + \beta \sum_{j \in N}
 p_{ij} v_{(1, j)}^*(\nu),
\]
which holds since $c_i \geq 0$ by Assumption \ref{ass:key} (recall that  $d_i \equiv 0$).
\end{proof}

We next give the main result of this section (whose proof
uses results from Section \ref{s:pclic}).
\begin{theorem}
\label{the:indxb}
The reformulated restless project 
is $\hat{\mathcal{F}}$-indexable, with MPI given by
\textup{(\ref{eq:nu1i})}. 
\end{theorem}
\begin{proof}
We first show that it is optimal to rest the project at $(1, i)$ iff
  $\nu \geq \nu_{(1, i)}^{\textup{AT}}$.
We have:
\begin{align*}
\text{it is optimal to rest  at $(1, i)$} &
\Longleftrightarrow
0 \geq \max_{S_0 \subseteq S_1 \subseteq N\colon i
 \in S_1} f^{S_0 \oplus S_1}_{(1, i)} - \nu 
 g^{S_0 \oplus S_1}_{(1, i)} \\
& \Longleftrightarrow \nu 
 \geq \max_{S_0 \subseteq S_1 \subseteq N\colon i
 \in S_1}
\frac{f^{S_0 \oplus S_1}_{(1, i)}}{g^{S_0 \oplus S_1}_{(1, i)}} \\
& \Longleftrightarrow \nu 
 \geq \max_{S_1 \subseteq N\colon i \in S_1} 
\frac{f_i^{S_1}}{g_i^{S_1}} =
\nu^{\textup{AT}}_{(1, i)},
\end{align*}
where we have used Lemma \ref{lma:intuit} and the following relations in Lemmas
\ref{lma:pwm}(b) and \ref{lma:hatrmp}(b) below: $g^{S_0 \oplus S_1}_{(1, i)} =
g_i^{S_1}$ and $f^{S_0 \oplus S_1}_{(1, i)} = f_i^{S_1}$ for 
 $S_0 \subseteq S_1 \subseteq N$ and $i
 \in S_1$. 

It remains to show that it is optimal to rest the project at $(0, i)$ iff 
  $\nu \geq \nu_{(0, i)}^{\textup{AT}}$.
We have:
\begin{align*}
\text{it is optimal to rest at $(0, i)$} &
\Longleftrightarrow
0 \geq \max_{S_0 \subseteq S_1 \subseteq N\colon i \in S_0} f^{S_0 \oplus S_1}_{(0, i)} - \nu 
 g^{S_0 \oplus S_1}_{(0, i)} \\
& \Longleftrightarrow
\nu \geq 
\max_{S_0 \subseteq S_1 \subseteq N\colon i \in S_0}
\frac{f^{S_0 \oplus S_1}_{(0, i)}}{g^{S_0 \oplus
    S_1}_{(0, i)}} \\
& \Longleftrightarrow
\nu \geq 
\max_{S_1 \subseteq N\colon i \in S_0}
\frac{f^{S_1}_{i} - c_i}{g^{S_1}_{i}} = 
\nu^{\textup{AT}}_{(0, i)},
\end{align*}
where we have used Lemma \ref{lma:intuit} and the following relations
in Lemmas \ref{lma:pwm}(c) and
\ref{lma:hatrmp}(c) below: $g^{S_0 \oplus S_1}_{(0, i)} =
g_i^{S_1}$ and $f^{S_0 \oplus S_1}_{(0, i)} = f_i^{S_1} - c_i$ for 
 $i \in S_0 \subseteq S_1 \subseteq N$. 

The above establishes that the project is indexable, with continuation
and switching MPIs given by
\textup{(\ref{eq:nu1i})}. The fact that it
is $\hat{\mathcal{F}}$-indexable then follows from 
Lemma \ref{lma:nu1gtnu0}.
\end{proof}

\subsection{On a Property of 
Optimal Policies}
\label{s:refatpop}
The main result of \citet{asatene96} establishes that their 
continuation and switching indices partially characterize
optimal policies for the MABPSC.

\begin{theorem}[\citet{asatene96}]
\label{the:at}
An optimal policy for the MABPSC has the property that 
decisions about which project to engage need to be made only at
stopping times that achieve the appropriate (continuation or
switching) index.
\end{theorem}

In \citet[Sec.\ 3.C]{asatene96} it is specified how
to construct the continuation (stopping) sets achieving the stopping times
in Theorem \ref{the:at}.
The next result revisits
 such an issue in light of the above MPI characterization of
the indices in (\ref{eq:nu1i}).
It follows immediately from the previous section's results, and hence
we omit the proof.

\begin{proposition}
\label{pro:opacsets} \mbox{ }
\begin{itemize}
\item[\textup{(a)}] The minimal continuation set
 achieving the continuation index 
$\nu_{(1, i)}^*$ in \textup{(\ref{eq:nu1i})} is
\[
S_1^*\left(\nu_{(1, i)}^*\right) = 
\big\{j \in N\colon \nu_{(1, j)}^* > \nu_{(1, i)}^*\big\}.
\]
\item[\textup{(b)}] The minimal continuation set
 achieving the switching index 
$\nu_{(0, i)}^*$ in \textup{(\ref{eq:nu1i})} is
\[
S_1^*\left(\nu_{(0, i)}^*\right) = 
\big\{j \in N\colon \nu_{(1, j)}^* > \nu_{(0, i)}^*\big\}.
\]
\end{itemize}
\end{proposition}

Notice that Proposition \ref{pro:opacsets} shows
that the MPI-based policy for the MABPSC, though generally not
optimal,
 satisfies the
property of optimal policies in
Theorem \ref{the:at}.

\subsection{PCL-indexability and MPI Algorithm}
\label{s:pclic}
We next outline the approach we will deploy to 
  compute the MPI,  based on the \emph{adaptive-greedy} MPI algorithm 
introduced in  \cite{nmaap01,nmmp02}.
Such an algorithm computes the MPI of restless bandits that
 satisfy a strong form of indexability termed 
\emph{PCL-indexability}.

Given  $a \in \{0, 1\}$ and 
$S_0 \oplus S_1 \in \hat{\mathcal{F}}$, denote by 
$\langle a, S_0 \oplus S_1\rangle$ the policy that takes action 
$a$ in the initial period and adopts the \emph{$S_0 \oplus S_1$-active
  policy} thereafter. 
For each augmented state $(a^-, i)$ and active set $S_0 \oplus S_1 \in
\hat{\mathcal{F}}$, define the \emph{marginal work measure}
\begin{equation}
\label{eq:hatismwm}
w^{S_0 \oplus S_1}_{(a^-, i)} 
\triangleq g^{\langle 1, S_0 \oplus S_1\rangle}_{(a^-, i)} - 
g^{\langle 0, S_0 \oplus S_1\rangle}_{(a^-, i)},
\end{equation}
along with  the \emph{marginal reward measure} 
\begin{equation}
\label{eq:hatismrm}
r^{S_0 \oplus S_1}_{(a^-, i)}
\triangleq f^{\langle 1, S_0 \oplus S_1\rangle}_{(a^-, i)} 
- f^{\langle 0,  S_0 \oplus S_1\rangle}_{(a^-, i)}
\end{equation}
and the
\emph{marginal productivity measure} 
\begin{equation}
\label{eq:hatismpm}
\nu^{S_0 \oplus S_1}_{(a^-, i)} \triangleq 
\frac{r^{S_0 \oplus S_1}_{(a^-, i)}}{w^{S_0 \oplus S_1}_{(a^-, i)}}.
\end{equation} 
We will see later (cf.\ Proposition \ref{pro:mwmpos}) that 
$w^{S_0 \oplus S_1}_{(a^-, i)} > 0$, ensuring that $\nu^{S_0 \oplus S_1}_{(a^-, i)}$ is well defined.

We will use the adaptive-greedy MPI algorithm
 $\mathrm{AG}_{\hat{\mathcal{F}}}$
 in Table \ref{fig:hatag2}, where
$n \triangleq |N|$ is the number of project states in the original
 formulation.
The algorithm's input 
consists of the project's parameters (which are
  not explicitly stated), along with the given active-set family 
$\hat{\mathcal{F}}$.
It  returns as output 
a string $\{(a_k^{-}, i_k)\}_{k=1}^{2 n}$ of 
distinct  augmented states spanning $\hat{N}$, which satisfy
$\hat{S}^k \triangleq \{(a_1^{-}, i_1), \ldots, (a_k^{-}, i_k)\}
\in \hat{\mathcal{F}}$ for $1 \leq k \leq 2 n$,
along with corresponding index values $\{\nu^*_{(a_k^{-},
  i_k)}\}_{k=1}^{2n}$.
Notice that each active set $\hat{S}^k$ is of the form 
$S_0 \oplus S_1$ for certain $S_0 \subseteq S_1 \subseteq N$. 
We use in this algorithm the compact form $\hat{S}^k$ for notational
convenience. 

\begin{table}[ht!]
\caption{Version 1 of MPI Algorithm $\mathrm{AG}_{\hat{\mathcal{F}}}$.}
\begin{center}
\fbox{%
\begin{minipage}{5in}
\textbf{ALGORITHM} $\mathrm{AG}_{\hat{\mathcal{F}}}$: \\
\textbf{Output:}
$\big\{(a_k^{-}, i_k), \nu^*_{(a_k^{-}, i_k)}\big\}_{k=1}^{2 n}$
\begin{tabbing}
$\hat{S}^0 := \emptyset$ \\
\textbf{for} \= $k := 1$ \textbf{to}  $2 n$ \textbf{do} \\
\> \textbf{pick}  
 $(a_k^-, i_k) \in \argmax
      \big\{\nu^{\hat{S}^{k-1}}_{(a^-, i)}\colon
                 (a^-, i) \in \hat{N} \setminus \hat{S}^{k-1}, 
   \hat{S}^{k-1} \cup \{(a^-, i)\} \in \hat{\mathcal{F}})\big\}$ \\
 \>  $\nu^*_{(a_k^-, i_k)} := 
 \nu^{\hat{S}^{k-1}}_{(a_k^-, i_k)}$;  \, $\hat{S}^{k} := \hat{S}^{k-1} \cup \big\{(a_k^-,
 i_k)\big\}$ \\
\textbf{end} \{ for \}
\end{tabbing}
\end{minipage}}
\end{center}
\label{fig:hatag2}
\end{table}

\begin{definition}[PCL($\hat{\mathcal{F}}$)-indexability]
\label{def:pclfi}
\textup{We say the project is
  \emph{PCL$(\hat{\mathcal{F}})$-indexable}  if: 
\begin{itemize}
\item[(i)] Positive marginal workloads:
  $w^{S_0 \oplus S_1}_{(a^-, i)} > 0$ for $(a^-, i) \in \hat{N}$
  and $S_0 \oplus S_1 \in \hat{\mathcal{F}}$.
\item[(ii)] Monotone nonincreasing index: algorithm $\mathrm{AG}_{\hat{\mathcal{F}}}$'s output satisfies
\[
\nu^*_{(a_1^-, i_1)} 
\geq \nu^*_{(a_2^-, i_2)} \geq \cdots
  \geq \nu^*_{(a_{2n}^-, i_{2n})}.
\]
\end{itemize}}
\end{definition}

The following result is established in  \citet[Cor.\ 2]{nmaap01}
and, in a more general setting whose formulation we adopt herein, in \citet[Th.\ 6.3]{nmmp02}.
\begin{theorem}
\label{the:sic}
A PCL$(\hat{\mathcal{F}})$-indexable restless project is
$\hat{\mathcal{F}}$-indexable, and the index $\nu^*_{(a^-, i)}$
computed by algorithm $\mathrm{AG}_{\hat{\mathcal{F}}}$ is its MPI.
\end{theorem}

While 
$\hat{\mathcal{F}}$-indexability was already established in Theorem
\ref{the:indxb} above,
we will invoke Theorem \ref{the:sic} in Section \ref{s:ppclfi} to ensure the validity
of algorithm $\mathrm{AG}_{\hat{\mathcal{F}}}$ to compute the MPI.
Notice that ties for picking the $(a_k^-,
i_k)$ in algorithm $\mathrm{AG}_{\hat{\mathcal{F}}}$ can be
broken arbitrarily.

\subsection{Version 2 of the MPI Algorithm}
\label{s:evaag}
The MPI algorithm 
in Table \ref{fig:hatag2} (Version 1) is readily reformulated into the more
explicit Version 2 in Table \ref{fig:expalg}, 
which uses the definition of 
$\hat{\mathcal{F}}$ in (\ref{eq:fhat}). 
We use in this and later versions a more convenient algorithm-like notation, writing,
e.g., 
$\nu^{S_0^{k_0-1} \oplus S_1^{k_1-1}}_{(0, j)}$ as $\nu^{(k_0-1, k_1-1)}_{(0, j)}$.
Notice that the active sets constructed in both versions along with 
respective counters $k$ and $k_0, k_1$  are related by
$\hat{S}^{k-1} = S_0^{k_0-1} \oplus S_1^{k_1-1}$
 with  $k = k_0 + k_1 - 1$ and $k_0 \leq k_1$.

Such a Version 2 decouples its output into the two augmented-state strings 
$\big\{(0, i_0^{k_0})\big\}_{k_0=1}^{n}$ and $\big\{(1, i_1^{k_1})\big\}_{k_1=1}^{n}$,
with
$N = \{i_0^{1}, \ldots, i_0^{n}\} = \{i_1^{1}, \ldots, i_1^{n}\}$,
along with corresponding switching and continuation index values.
The active sets $S_0^{k_0}$ and $S_1^{k_1}$ in the
algorithm are given by
$S_0^{k_0} = \{i_0^{1}, \ldots, i_0^{k_0}\}$ 
and 
$S_1^{k_1} = \{i_1^{1}, \ldots, i_1^{k_1}\}$, 
and satisfy 
$S_0^{k_0} \subset S_1^{k_1}$ for  $1 \leq k_0 < k_1 \leq n$.

We will later show  that Version 2
can still be substantially simplified,  decoupling 
the computation of the continuation and switching MPIs into two
 algorithms. 
Yet to achieve such simplifications we will need the
PCL-indexability analysis  carried out in the next section.

\begin{table}
\caption{Version 2 of MPI Algorithm
$\mathrm{AG}_{\hat{\mathcal{F}}}$.}
\begin{center}
\fbox{%
\begin{minipage}{5in}
\textbf{ALGORITHM} $\mathrm{AG}_{\hat{\mathcal{F}}}$: \\
\textbf{Output:}
$\big\{(0, i_0^{k_0}), \nu^*_{(0, i_0^{k_0})}\big\}_{k_0=1}^{n}$, 
$\big\{(1, i_1^{k_1}), \nu^*_{(1, i_1^{k_1})}\big\}_{k_1=1}^{n}$
\begin{tabbing}
$S_0^{0} := \emptyset$; \, $S_1^{0} := \emptyset$; \,
$k_0 := 1$; \, $k_1 := 1$ \\
\textbf{while} \= $k_0+k_1 \leq 2 n+1$
 \textbf{do} \\
\> \textbf{if } $k_1 \leq n$ \textbf{pick} 
 $j_1^{\max} \in \argmax \big\{\nu^{(k_0-1, k_1-1)}_{(1, j)}\colon j \in N \setminus S_1^{k_1-1}\big\}$ \\
\> \textbf{if } $k_0 < k_1$ \textbf{pick} 
 $j_0^{\max} \in \argmax \big\{\nu^{(k_0-1, k_1-1)}_{(0, j)}\colon j \in  S_1^{k_1-1} \setminus S_0^{k_0-1}\big\}$ \\
\> \textbf{if } \= $k_1 = n+1$ \textbf{or} $\big\{k_0 < k_1 \leq n
\textbf{ and } \nu^{(k_0-1, k_1-1)}_{(1, j_1^{\max})} < 
  \nu^{(k_0-1, k_1-1)}_{(0, j_0^{\max})}\big\}$ \\
\> \> 
   $i_0^{k_0} := j_0^{\max}$; \,  
$\nu^*_{(0, i_1^{k_0})} := \nu^{(k_0-1, k_1-1)}_{(0,
  i_1^{k_0})}$; \, $S_0^{k_0} := S_0^{k_0-1} \cup \{i_0^{k_0}\}$; \,  
  $k_0 := k_0 + 1$ \\
\> \textbf{else}           \\
\> \> 
$i_1^{k_1} := j_1^{\max}$; \, 
$\nu^*_{(1, i_1^{k_1})} := 
\nu^{(k_0-1, k_1-1)}_{(1, i_1^{k_1})}$; \, 
$S_1^{k_1} := S_1^{k_1-1} \cup \{i_1^{k_1}\}$;  \,
  $k_1 := k_1 + 1$ \\
\> \textbf{end }  \{ if \} \\
\textbf{end }  \{ while \} 
\end{tabbing}
\end{minipage}}
\end{center}
\label{fig:expalg}
\end{table}

\section{PCL-Indexability Analysis}
\label{s:pclapsc}
We set out in this section to carry out a PCL-indexability analysis 
of a single project with switching costs as above, in its restless reformulation.

\subsection{Work and Marginal Work Measures}
\label{s:wmwm}
We start by addressing calculation of work and marginal work 
measures $g^{S_0 \oplus
  S_1}_{(a^-, i)}$ and $w^{S_0 \oplus S_1}_{(a^-, i)}$ for $(a^-, i) \in \hat{N}$ and $S_0 \oplus S_1
  \in \hat{\mathcal{F}}$.
We will show that they are closely related to corresponding
   measures $g^S_i$ and $w^S_i$ for the 
 underlying nonrestless project, where
stationary deterministic policies are represented
by their active sets $S \subseteq N$.

For each $S \subseteq N$, work measures $g_i^S$ are
characterized by the evaluation equations

\begin{equation}
\label{eq:eegis}
g^S_i = 
\begin{cases}
1 + \beta \sum_{j \in S} p_{i j} g^S_j & \text{if } i \in S \\
0 & \text{otherwise}.
\end{cases}
\end{equation}
Note that (\ref{eq:eegis})  can
be reduced to a linear system with coefficient matrix 
$\mathbf{I}_S - \beta \mathbf{P}_{SS}$, where $\mathbf{I}_S$ is the
identity matrix indexed by $S$ and $\mathbf{P}_{SS} \triangleq
(p_{ij})_{i, j \in S}$. Its solution is unique, as matrix theory ensures
that a matrix of the form $\mathbf{I}_S - \beta \mathbf{P}_{SS}$ is
invertible if $\mathbf{P}_{SS}$ is substochastic (i.e., it has
nonnegative elements and its rows add up to at most
unity) and $0 < \beta < 1$.

Further, the marginal work measure $w^S_i$ is evaluated by 
\begin{equation}
\label{eq:ismwm}
\begin{split}
w^S_i & \triangleq g^{\langle 1, S\rangle}_i - g^{\langle 0,
  S\rangle}_i = 1 + \beta \sum_{j \in N} p_{ij} g^S_j - \beta g^S_i = 
\begin{cases}
(1-\beta) g^S_i & \text{if } i \in S \\
1 + \beta \sum_{j \in S} p_{ij} g^S_j & \textrm{otherwise.}
\end{cases}
\end{split}
\end{equation}
Notice that (\ref{eq:ismwm}) implies
\begin{equation}
\label{eq:wsipos}
w^S_i > 0, \quad i \in N.
\end{equation}

We now return to the project's
 restless reformulation.  The following result
gives the evaluation equations for work measure $g^{S_0 \oplus
  S_1}_{(a^-, i)}$, for a given active set
$S_0 \oplus S_1 \in \hat{\mathcal{F}}$.
\begin{lemma}
\label{lma:eewm}
\begin{equation*}
g^{S_0 \oplus S_1}_{(a^-, i)} = 
\begin{cases}
1 + \beta \sum_{j \in N} p_{ij} g^{S_0 \oplus S_1}_{(1, j)} & \textup{if } \, i \in S_{a^-} \\ 
\beta  g^{S_0 \oplus S_1}_{(0, i)} & \textup{otherwise}.
\end{cases}
\end{equation*}
\end{lemma}

The next result represents work measure
$g^{S_0 \oplus S_1}_{(a^-, i)}$ in terms of the 
$g^S_i$.
\begin{lemma}
\label{lma:pwm}
For $S_0 \oplus S_1 \in \hat{\mathcal{F}}$\textup{:}
\begin{itemize}
\item[\textup{(a)}] $g^{S_0 \oplus S_1}_{(a^-, i)} = g^{S_1}_i = 0$, for
  $a^- \in \{0, 1\}, i \in N \setminus S_1$.
\item[\textup{(b)}] $g^{S_0 \oplus S_1}_{(1, i)} = 
g^{S_1}_i$, for $i \in S_1$.
\item[\textup{(c)}] $g^{S_0 \oplus S_1}_{(0, i)} = 
g^{S_1}_i$, for $i \in S_0$.
\item[\textup{(d)}] $g^{S_0 \oplus S_1}_{(0, i)} = 0$, 
for $i \in S_1 \setminus S_0$.
\end{itemize}
\end{lemma}
\begin{proof}
(a) This part follows immediately from the definition of policy 
$S_0 \oplus S_1$.

(b) For $i \in S_1$, we can write
\begin{align*}
g^{S_0 \oplus S_1}_{(1, i)} & = 1 + \beta \sum_{j \in S_1} p_{ij}
g^{S_0 \oplus S_1}_{(1, j)}+ \beta \sum_{j \in N \setminus S_1} 
p_{ij}
g^{S_0 \oplus S_1}_{(1, j)} = 1 + \beta \sum_{j \in S_1} p_{ij} g^{S_0 \oplus S_1}_{(1, j)},
\end{align*}
where we have used Lemma \ref{lma:eewm} and part (a).
Hence, the $g^{S_0 \oplus S_1}_{(1, i)}$ satisfy the evaluation
equations in (\ref{eq:eegis}) for the $g^{S_1}_i$,  for $i
\in S_1$, which yields the result.

(c) For $i \in S_0$,  we have
\begin{align*}
g^{S_0 \oplus S_1}_{(0, i)} & = 1 +  \beta \sum_{j \in S_1} p({i,
  j})
g^{S_0 \oplus S_1}_{(1, j)} + \beta \sum_{j \in N \setminus S_1} p_{ij}
g^{S_0 \oplus S_1}_{(1, j)} = g^{S_0 \oplus S_1}_{(1, i)} = g^{S_1}_i,
\end{align*}
where we have used Lemma \ref{lma:eewm}, part (b) and the relation $S_0 \subseteq S_1$.

(d) This part follows immediately from the definition of policy 
$S_0 \oplus S_1$.
\end{proof}

Regarding  $w^{S_0 \oplus S_1}_{(a^-, i)}$, we draw on 
(\ref{eq:hatismwm}) and Lemma \ref{lma:eewm} to obtain
\begin{equation}
\label{eq:ws0s1oi}
w^{S_0 \oplus S_1}_{(0, i)} =
w^{S_0 \oplus S_1}_{(1, i)} = 
1 + \beta \sum_{j \in N} p_{ij} g^{S_0 \oplus S_1}_{(1, j)} - 
\beta  g^{S_0 \oplus S_1}_{(0, i)}.
\end{equation}

The following result represents 
$w^{S_0 \oplus S_1}_{(a^-, i)}$ in terms of the $w^S_i$ in
(\ref{eq:ismwm}).

\begin{lemma}
\label{lma:hatwis}
For $a^- \in \{0, 1\}, S_0 \oplus S_1 \in \hat{\mathcal{F}}$\textup{:}
\begin{itemize}
\item[\textup{(a)}] $w^{S_0 \oplus S_1}_{(a^-, i)} = 
  w^{S_1}_i$,  for $i \in S_0 \cup N \setminus S_1$.
\item[\textup{(b)}] $w^{S_0 \oplus S_1}_{(a^-, i)} = 
  w^{S_1}_i/(1-\beta)$,  for $i \in S_1
  \setminus S_0$.
\end{itemize}
\end{lemma}
\begin{proof}
(a) We can write, for $i \in S_0 \cup N \setminus S_1$,
\begin{align*}
w^{S_0 \oplus S_1}_{(a^-, i)} & = 
1 + \beta \sum_{j \in N} p_{ij} g^{S_0 \oplus S_1}(1, j)  - 
\beta g^{S_0 \oplus S_1}_{(0, i)}   
 = 1 + \beta \sum_{j \in S_1} p_{ij} g^{S_1}_j - \beta g^{S_1}_i 
 = w^{S_1}_i,
\end{align*}
where we have used (\ref{eq:ws0s1oi}), 
Lemma \ref{lma:pwm}(a, b, c) and (\ref{eq:ismwm}).

(b) We can write, for $i \in S_1 \setminus S_0$,
\begin{align*}
w^{S_0 \oplus S_1}_{(a^-, i)} & = 
1 + \beta \sum_{j \in N} p_{ij} g^{S_0 \oplus S_1}(1, j)  - 
\beta g^{S_0 \oplus S_1}_{(0, i)}   
 = 1 + \beta \sum_{j \in S_1} p_{ij} g^{S_1}_j 
 = w^{S_1}_i + \beta g^{S_1}_i  = w^{S_1}_i/(1-\beta),
\end{align*}
where we have used (\ref{eq:ws0s1oi}), 
Lemma \ref{lma:pwm}(a, b, d) and (\ref{eq:ismwm}).
\end{proof}

From the above, we obtain the following result (cf.\
Definition \ref{def:pclfi}(i)).
\begin{proposition}
\label{pro:mwmpos}
$w^{S_0 \oplus S_1}_{(a^-, i)} > 0$, for $(a^-, i) \in \hat{N}, S_0 \oplus S_1 \in \hat{\mathcal{F}}$.
\end{proposition}
\begin{proof}
The result follows immediately from 
(\ref{eq:wsipos}) via  Lemma \ref{lma:hatwis}.
\end{proof}

\subsection{Reward and Marginal Reward Measures}
\label{s:rmwm}
We next address the calculation of 
measures $f^{S_0 \oplus
  S_1}_{(a^-, i)}$ and $r^{S_0 \oplus S_1}_{(a^-, i)}$, relating
them to their counterparts
    $f^S_i$ and $r^S_i$
  for the underlying nonrestless
project with no
startup costs.

For each active set $S \subseteq N$, 
the reward measure $f^S_i$
is characterized by the equations

\begin{equation}
\label{eq:nteefis}
f^S_i =
\begin{cases}
R_i + \beta \sum_{j \in S} p_{ij} f^S_j & \text{if } i \in S \\
0 & \text{otherwise},
\end{cases}
\end{equation}
whose solution is unique, while the
marginal reward measure $r^S_i$ is given by
\begin{equation}
\label{eq:ntsimpris}
\begin{split}
r^S_i & \triangleq f^{\langle 1, S\rangle}_i - f^{\langle 0,
  S\rangle}_i 
= 
R_i + \beta \sum_{j \in S} p_{ij} f^S_j - \beta f^S_i
   =
\begin{cases}
(1-\beta) f^S_i  & \text{if } i \in S \\ 
R_i + \beta \sum_{j \in S} p_{ij} f^S_j  & \text{otherwise.}
\end{cases}
\end{split}
\end{equation}

Returning to the  restless formulation, the next result
gives the evaluation equations for reward measures $f^{S_0 \oplus
  S_1}_{(a^-, i)}$, for a given active set
$S_0 \oplus S_1 \in \hat{\mathcal{F}}$. Recall the notation in
(\ref{eq:amminus}).

\begin{lemma}
\label{lma:eefm}
\begin{equation*}
f^{S_0 \oplus S_1}_{(a^-, i)} = 
\begin{cases}
R_i -  (1-a^-) c_i + \beta \sum_{j \in N} p_{ij} f^{S_0 \oplus S_1}_{(1, j)} & \textup{if } i \in S_{a^-} \\ 
\beta  f^{S_0 \oplus S_1}_{(0, i)} & \textup{otherwise}.
\end{cases}
\end{equation*}
\end{lemma}

The next result represents reward measure
$f^{S_0 \oplus S_1}_{(a^-, i)}$ in terms of the 
$f^S_i$.

\begin{lemma}
\label{lma:hatrmp}
For $S_0 \oplus S_1 \in \hat{\mathcal{F}}$\textup{:}
\begin{itemize}
\item[\textup{(a)}] $f^{S_0 \oplus S_1}_{(a^-, i)} = 0 = f^{S_1}_i$, 
for $a^- \in \{0, 1\}, i \in N \setminus S_1$.
\item[\textup{(b)}] $f^{S_0 \oplus S_1}_{(1, i)} = f^{S_1}_i$, 
for $i \in S_1$.
\item[\textup{(c)}] $f^{S_0 \oplus S_1}_{(0, i)} = f^{S_1}_i - c_i$, 
for $i \in S_0$.
\item[\textup{(d)}] $f^{S_0 \oplus S_1}_{(0, i)} = 0 = 
f^{S_0}_i$, 
for $i \in S_1 \setminus S_0$.
\end{itemize}
\end{lemma}
\begin{proof}
(a) This part follows immediately from the definition of policy 
$S_0 \oplus S_1$.

(b) We can write, for $i \in S_1$,
\begin{align*}
f^{S_0 \oplus S_1}_{(1, i)} & =
R_i + \beta \sum_{j \in S_1} p_{ij} f^{S_0 \oplus S_1}_{(1, j)}
 + \beta \sum_{j \in N \setminus S_1} p_{ij} f^{S_0 \oplus S_1}_{(1,
 j)} = R_i + \beta \sum_{j \in S_1} p_{ij} f^{S_0 \oplus S_1}_{(1, j)},
\end{align*}
where we have used Lemma \ref{lma:eefm} and part (a). Hence, the 
$f^{S_0 \oplus S_1}_{(1, i)}$, for $i \in S_1$, satisfy the
evaluation equations in (\ref{eq:nteefis}) for corresponding
terms $f^{S_1}_i$, which yields the result.

(c) We can write, for $i \in S_0$,
\begin{align*}
f^{S_0 \oplus S_1}_{(0, i)} & = 
R_i - c_i + \beta \sum_{j \in S_1} p_{ij} f^{S_0 \oplus S_1}_{(1, j)}
 + \beta \sum_{j \in N \setminus S_1} p_{ij} f^{S_0 \oplus S_1}(1,
 j) 
 = f^{S_0 \oplus S_1}_{(1, i)} - c_i = f^{S_1}_i  - c_i,
\end{align*}
where we have used  $S_0 \subseteq S_1$ along with parts (a, b).

(d) This part follows immediately from the definition of policy 
$S_0 \oplus S_1$.
\end{proof}

Regarding marginal reward measure $r^{S_0 \oplus S_1}_{(a^-, i)}$, 
we draw on (\ref{eq:hatismrm}) and Lemma \ref{lma:eefm} to obtain
\begin{equation}
\label{eq:rs0s1oi}
r^{S_0 \oplus S_1}_{(a^-, i)} = 
R_i -  (1-a^-) c_i + \beta \sum_{j \in N} p_{ij} f^{S_0
  \oplus S_1}_{(1, j)} -
\beta  f^{S_0 \oplus S_1}_{(0, i)}.
\end{equation}

The following result represents
$r^{S_0 \oplus S_1}_{(a^-, i)}$ in terms of the $r^S_i$.

\begin{lemma}
\label{lma:pmrm}
For $S_0 \oplus S_1 \in \hat{\mathcal{F}}$\textup{:}
\begin{itemize}
\item[\textup{(a)}] 
$r^{S_0 \oplus S_1}_{(0, i)} = r^{S_0 \oplus S_1}_{(1, i)} - c_i$, for 
 $i \in N$.
\item[\textup{(b)}] $r^{S_0 \oplus S_1}_{(1, i)} =
r^{S_1}_i$,  for $i \in N \setminus S_1$.
\item[\textup{(c)}] $r^{S_0 \oplus S_1}_{(1, i)} =
 r^{S_1}_i + \beta c_i$, for $i \in S_0$.
\item[\textup{(d)}] $r^{S_0 \oplus S_1}_{(1, i)} =
r^{S_1}_i/(1-\beta)$, for $i \in S_1 \setminus S_0$.
\end{itemize}
\end{lemma}
\begin{proof}
(a) This part follows immediately from (\ref{eq:rs0s1oi}).

(b) We can write, for $i \in N \setminus S_1$,
\begin{align*}
r^{S_0 \oplus S_1}_{(1, i)} & = 
R_i + \beta \sum_{j \in N} p_{ij} f^{S_0 \oplus S_1}_{(1, j)}
-  f^{S_0 \oplus S_1}_{(1, i)} 
= R_i + \beta \sum_{j \in S_1} p_{ij} f^{S_1}_j - f^{S_1}_i 
 = r^{S_1}_i,
\end{align*}
where we have used (\ref{eq:rs0s1oi}), Lemma \ref{lma:hatrmp}(a, b),
and (\ref{eq:ntsimpris}).

(c) We can write, for $i \in S_0$, 
\[
r^{S_0 \oplus S_1}_{(1, i)} = 
R_i + \beta \sum_{j \in N} p_{ij} f^{S_0 \oplus S_1}_{(1, j)}
-  \beta f^{S_0 \oplus S_1}_{(0, i)} 
= R_i + \beta \sum_{j \in S_1} p_{ij} f^{S_1}_j - \beta (f^{S_1}_i - c_i)
 = r^{S_1}_i + \beta c_i,
\]
where we have used $S_0 \subseteq S_1$, (\ref{eq:rs0s1oi}), 
Lemma \ref{lma:hatrmp}(a, b, c) and (\ref{eq:ntsimpris}).

(d) We can write, for $i \in S_1 \setminus S_0$,
\begin{align*}
r^{S_0 \oplus S_1}_{(1, i)} & = 
R_i + \beta \sum_{j \in N} p_{ij} f^{S_0 \oplus S_1}_{(1, j)} - \beta f^{S_0 \oplus S_1}_{(0, i)} = f_i^{S_1} = \frac{r^{S_1}_i}{1-\beta},
\end{align*}
where we have used (\ref{eq:rs0s1oi}), 
Lemma \ref{lma:hatrmp}(a, b, d),
(\ref{eq:nteefis}) and (\ref{eq:ntsimpris}). 
This completes the proof.
\end{proof}

\subsection{Marginal Productivity Measures}
\label{s:rmpm}
We continue to address calculation of marginal productivity
measures $\nu^{S_0 \oplus
  S_1}_{(a^-, i)}$  in (\ref{eq:hatismpm}), 
relating them to their
counterparts
  for the underlying nonrestless project, given by 
\begin{equation}
\label{eq:nusi}
\nu^S_i \triangleq \frac{r^S_i}{w^S_i}, \quad i \in N, S \subseteq N.
\end{equation}

The next result represents the required
$\nu^{S_0 \oplus S_1}_{(a^-, i)}$ in terms of the 
$\nu^{S}_i$.

\begin{lemma}
\label{lma:pmpr}
For $S_0 \oplus S_1 \in \hat{\mathcal{F}}$\textup{:}
\begin{itemize}
\item[\textup{(a)}] 
$\nu^{S_0 \oplus S_1}_{(0, i)} = \nu^{S_0 \oplus S_1}_{(1, i)} - c_i/w^{S_0 \oplus S_1}_{(1, i)}$, for 
 $i \in N$.
\item[\textup{(b)}] $\nu^{S_0 \oplus S_1}_{(1, i)} =
\nu^{S_1}_i = \nu^{\emptyset \oplus S_1}_{(1, i)}$,  for $i \in N
\setminus S_0$.
\item[\textup{(c)}] $\nu^{S_0 \oplus S_1}_{(1, i)} =
 \nu^{S_1}_i + \beta c_i / w^{S_1}_i$, for $i \in S_0$.
\item[\textup{(d)}] $\nu^{S_0 \oplus S_1}_{(0, i)} = 
\nu^{S_1}_i - (1-\beta) c_i
 / w^{S_1}_i = \nu^{\emptyset \oplus S_1}_{(0, i)}$, $i \in S_1 \setminus S_0$.
\end{itemize}
\end{lemma}
\begin{proof}
All parts follow immediately from (\ref{eq:hatismpm}),
(\ref{eq:nusi}), Lemma \ref{lma:hatwis} and Lemma \ref{lma:pmrm}.
\end{proof}

The above results allow us to
 reformulate Version 2 of the MPI algorithm  into 
the even more explicit
 Version 3 shown in Table \ref{fig:v3mpia}.
As before, we write, e.g., $\nu^{S_{k_0}^0 \oplus S_{k_1}^1}_{(0, j)}$ as
$\nu^{(k_0, k_1)}_{(0, j)}$, and $\nu^{S_{k_1}}_j$ as
 $\nu^{(k_1)}_j$.
Notice that in Version 3 we use $\nu^{(0, k_1-1)}_{(0, j)}$ (which denotes $\nu_{(0, j)}^{S_0
  \oplus S_{k_1-1}}$) in place of $\nu^{(k_0-1, k_1-1)}_{(0, j)}$, 
drawing on Lemma \ref{lma:pmpr}(d). We do so for
computational reasons, as storage of quantities 
$\nu^{(0,  k_1-1)}_{(0, j)}$ requires one less dimension than storage of the 
$\nu^{(k_0-1, k_1-1)}_{(0, j)}$.

\begin{table}[ht!]
\caption{Version 3 of Algorithm $\mathrm{AG}_{\hat{\mathcal{F}}}$.}
\begin{center}
\fbox{%
\begin{minipage}{5in}
\textbf{ALGORITHM} $\mathrm{AG}_{\hat{\mathcal{F}}}$: \\
\textbf{Output:}
$\big\{(0, i_0^{k_0}), \nu^*_{(0, i_0^{k_0})}\big\}_{k_0=1}^{n}$, 
$\big\{i_1^{k_1}, \nu^*_{i_1^{k_1}}\big\}_{k_1=1}^{n}$
\begin{tabbing}
$S_0^{0} := \emptyset$; \, $S_1^{0} := \emptyset$; \, 
$k_0 := 1$; \, $k_1 := 1$ \\
\textbf{while} \= $k_0+k_1 \leq 2 n+1$
 \textbf{do} \\
\> \textbf{if } $k_1 \leq n$ \textbf{pick} 
 $j_1^{\max} \in \argmax \big\{\nu^{(k_1-1)}_j\colon j \in N \setminus S_1^{k_1-1}\big\}$ \\
\> $\nu^{(0, k_1-1)}_{(0, j)} := \nu^{(k_1-1)}_j 
  - (1-\beta) c_j / w^{(k_1-1)}_j, \, j \in  S_1^{k_1-1} \setminus S_0^{k_0-1}$ \\
\> \textbf{if } \= $k_0 < k_1$  \textbf{pick} 
 $j_0^{\max} \in \argmax \big\{\nu^{(0, k_1-1)}_{(0, j)}\colon j \in  S_1^{k_1-1} \setminus S_0^{k_0-1}\big\}$  \\
\> \textbf{if } \= $k_1 = n+1$ \textbf{or} $\big\{k_0 < k_1 \leq n
\textbf{ and } \nu^{(k_1-1)}_{j_1^{\max}} < \nu^{(0, k_1-1)}_{(0, j_0^{\max})}\big\}$ \\
\> \>  \= $i_0^{k_0} := j_0^{\max}$; \, 
 $\nu^*_{(0, i_0^{k_0})} := \nu^{(0, k_1-1)}_{(0, i_0^{k_0})}$; \,
$S_0^{k_0} := S_0^{k_0-1} \cup \{i_0^{k_0}\}$; \,  
  $k_0 := k_0 + 1$ \\
\> \textbf{else}   
 \\
\> \> \> $i_1^{k_1} := j_1^{\max}$; \, $\nu^*_{i_1^{k_1}} :=
\nu^{(k_1-1)}_{i_1^{k_1}}$; \,
$S_1^{k_1} := S_1^{k_1-1} \cup \{i_1^{k_1}\}$; \, 
  $k_1 := k_1 + 1$ \\
\> \textbf{end }  \{ if \} \\
\textbf{end }  \{ while \} 
\end{tabbing}
\end{minipage}}
\end{center}
\label{fig:v3mpia}
\end{table}

\subsection{Proving $PCL(\hat{\mathcal{F}})$-indexability}
\label{s:ppclfi}
We are now ready to establish $PCL(\hat{\mathcal{F}})$-indexability (cf.\
Theorem \ref{the:sic}).
\begin{theorem}
\label{the:pclfi}
The reformulated restless project  is
PCL$(\hat{\mathcal{F}})$-indexable.
\end{theorem}
\begin{proof}
Definition \ref{def:pclfi}(i) is established in
Proposition \ref{pro:mwmpos}. It remains to prove condition (ii),
  stating that the successive index values computed by algorithm
  $\mathrm{AG}_{\hat{\mathcal{F}}}$ are nonincreasing. 
Referring to Version 1 of the algorithm,
we will thus show that the $k$th and $(k+1)$th computed index values satisfy 
$\nu_{(a_k^-, i_k)}^* \geq \nu_{(a_{k+1}^-, i_{k+1})}^*$ for $1
  \leq k < 2n$.
Now, letting $\hat{S}^{k-1}$ and $\hat{S}^k =
  \hat{S}^{k-1} \cup \{(a_k^-, i_k)\}$ be as in Table \ref{fig:hatag2},
we use \citet[Prop.\ 6.4(b, c)]{nmmp02} to write
\[
r_{(a^-, i)}^{\hat{S}^k} - r_{(a^-, i)}^{\hat{S}^{k-1}} =
\frac{r_{(a_k^-, i_k)}^{\hat{S}^{k-1}}}{w_{(a_k^-,
    i_k)}^{\hat{S}^{k-1}}} \big(w_{(a^-, i)}^{\hat{S}^k} -
w_{(a^-, i)}^{\hat{S}^{k-1}}\big), \quad (a^-, i) \in \hat{N},
\]
which is immediately reformulated using (\ref{eq:hatismpm}) and
$\nu_{(a_k^-, i_k)}^* = \nu_{(a_k^-, i_k)}^{\hat{S}^{k-1}}$ as
\[
\nu_{(a^-, i)}^{\hat{S}^k} = \nu_{(a_k^-, i_k)}^* - 
\frac{w_{(a^-, i)}^{\hat{S}^{k-1}}}{w_{(a^-,
    i)}^{\hat{S}^{k}}} \big(\nu_{(a_k^-, i_k)}^* - \nu_{(a^-,
  i)}^{\hat{S}^{k-1}}\big), \quad (a^-, i) \in \hat{N}.
\]
Now, taking $(a^-, i) = (a_{k+1}^-, i_{k+1})$ in the latter identity
and using Proposition \ref{pro:mwmpos}
we obtain
\[
\nu_{(a_k^-, i_k)}^* \geq \nu_{(a_{k+1}^-, i_{k+1})}^*
\Longleftrightarrow 
\nu_{(a_k^-, i_k)}^* \geq \nu_{(a_{k+1}^-,
  i_{k+1})}^{\hat{S}^{k-1}},
\]
which shows that it suffices to prove that $\nu_{(a_k^-, i_k)}^* \geq 
\nu_{(a_{k+1}^-,  i_{k+1})}^{\hat{S}^{k-1}}$.

For such a purpose, recall that the algorithm picks $(a_{k}^-, i_{k})$ and $(a_{k+1}^-, i_{k+1})$ so that 
\begin{align*}
(a_k^-, i_k) & \in \argmax
      \big\{\nu^{\hat{S}^{k-1}}_{(a^-, i)}\colon
                 (a^-, i) \in N \setminus \hat{S}^{k-1}, 
   \hat{S}^{k-1} \cup \{(a^-, i)\} \in \hat{\mathcal{F}})\big\}, \\
(a_{k+1}^-, i_{k+1}) & \in \argmax
      \big\{\nu^{\hat{S}^{k}}_{(a^-, i)}\colon
                 (a^-, i) \in N \setminus \hat{S}^{k}, 
   \hat{S}^{k} \cup \{(a^-, i)\} \in \hat{\mathcal{F}})\big\},
\end{align*}
and hence there are two cases to consider. If
$\hat{S}^{k-1} \cup \{(a_{k+1}^-, i_{k+1})\} \in \hat{\mathcal{F}}$, the above relations yield $\nu_{(a_k^-,
  i_k)}^* = \nu_{(a_k^-, i_k)}^{\hat{S}^{k-1}} \geq \nu_{(a_{k+1}^-,
  i_{k+1})}^{\hat{S}^{k-1}}$, as required.
If $\hat{S}^{k-1} \cup \{(a_{k+1}^-,
i_{k+1})\} \not\in \hat{\mathcal{F}}$, since $\hat{S}^{k} \cup \{(a_{k+1}^-,
i_{k+1})\} \in \hat{\mathcal{F}}$, the structure of
$\hat{\mathcal{F}}$ as defined in (\ref{eq:fhat}) 
implies that it must be $a_k^- = 1$ and $(a_{k+1}^-, i_{k+1}) = (0,
i_k)$. 
We now use Lemma \ref{lma:pmpr}(a) along with the latter
identities to write
\[
\nu_{(a_k^-, i_k)}^* = \nu_{(1, i_k)}^* = \nu_{(1, i_k)}^{\hat{S}^{k-1}} = \nu_{(0, i_k)}^{\hat{S}^{k-1}} + c_i / w_{(1, i_k)}^{\hat{S}^{k-1}} \geq 
\nu_{(0,  i_k)}^{\hat{S}^{k-1}} = \nu_{(a_{k+1}^-,  i_{k+1})}^{\hat{S}^{k-1}},
\]
where the inequality follows from Assumption \ref{ass:key} and
Proposition \ref{pro:mwmpos}.
\end{proof}

\section{Decoupled Computation of the MPI}
\label{s:ecsi}
We set out in this section to further simplify Version 3 of the 
MPI algorithm
$\mathrm{AG}_{\hat{\mathcal{F}}}$ by decoupling 
computation of the continuation and the switching MPIs in a two-stage scheme.

\subsection{First Stage: Computing the Continuation Index}
\label{s:cmpi}
We start with the continuation index $\nu^*_{(1, i)}$, which is (cf.\ 
Section \ref{s:dpaidmpi}) the Gittins index
$\nu^*_i$
of the project without switching costs.
We will need additional quantities to feed
the second-stage switching-index algorithm discussed in Section
\ref{s:smpi}.

\begin{table}[ht!]
\caption{Gittins-Index Algorithmic Scheme $\mathrm{AG}^1$.}
\begin{center}
\fbox{%
\begin{minipage}{4in}
\textbf{ALGORITHM} $\mathrm{AG}^1$: \\
\textbf{Output:}
 $\big\{i_1^{k_1}, \nu^*_{i_1^{k_1}}, (w^{(k_1)}_j, \nu^{(k_1)}_j),  j \in S_1^{k_1}\big\}_{k_1=1}^n$

\begin{tabbing}
\textbf{set }
$S_1^{0} := \emptyset$; \textbf{compute} 
$\{(w^{(0)}_i, \nu^{(0)}_i)\colon i \in N\}$ \\
\textbf{for} \= $k_1 := 1$ \textbf{to} $n$ \textbf{do} \\
\> \textbf{pick} 
 $i_1^{k_1} \in \argmax
      \big\{\nu^{(k_1-1)}_i\colon
                i \in N \setminus S_1^{k_1-1}\big\}$ \\
 \> $\nu^*_{i_1^{k_1}} := 
 \nu^{(k_1-1)}_{i_1^{k_1}}$; \, $S_1^{k_1} := S_1^{k_1-1}
 \cup \{i_1^{k_1}\}$ \\
\> \textbf{compute} 
$\big\{(w^{(k_1)}_i, \nu^{(k_1)}_i)\colon i \in N\big\}$ \\
\textbf{end} 
\end{tabbing}
\end{minipage}}
\end{center}
\label{fig:ag2}
\end{table}

To compute the Gittins index and additional quantities we use the
 algorithmic scheme $\mathrm{AG}^1$ in Table \ref{fig:ag2}. 
This is a variant of the algorithm of 
\cite{vawabu}, reformulated as in \cite{nmijoc06}. See the latter paper for implementations such as algorithm FP(1),  which
performs $(4/3) n^3 + O(n^2)$ arithmetic operations, or the
 \emph{complete-pivoting algorithm} CP --- which we have found to be
 faster in practice despite its worse operation count

\subsection{Second Stage: Computing the Switching Index}
\label{s:smpi}
We next address computation of the switching index.
Consider algorithm $\mathrm{AG}^0$ in Table \ref{fig:hatag3}, which is fed the output of 
 $\mathrm{AG}^1$ and produces a sequence of states $i_0^{k_0}$ spanning
$N$, along with nonincreasing 
index values $\nu^*_{(0, i_0^{k_0})}$.
We can now give  the main result of this paper.
\begin{theorem}
\label{the:mpi0a}
Algorithm $\mathrm{AG}^0$ computes the switching index $\nu^*_{(0, i)}$.
\end{theorem}
\begin{proof}
The result follows by noticing that algorithm
$\mathrm{AG}^0$ is obtained from Version 3 of 
MPI algorithm $\mathrm{AG}_{\hat{\mathcal{F}}}$ in Table \ref{fig:v3mpia} by 
decoupling the computation of the $\nu^*_{(0, i)}$ and 
the $\nu_i^*$.
\end{proof}

\begin{table}[ht!]
\caption{Switching-Index Algorithm $\mathrm{AG}^0$.}
\begin{center}
\fbox{%
\begin{minipage}{4.5in}
\begin{tabbing}
\textbf{ALGORITHM} $\mathrm{AG}^0$: \\
\textbf{Input:} \= 
 $\big\{i_1^{k_1}, \nu^*_{i_1^{k_1}}, (w^{(k_1)}_j, \nu^{(k_1)}_j),  j \in S_1^{k_1}\big\}_{k_1=1}^n$ \\
\textbf{Output:}
$\big\{i_0^{k_0}, \nu^*_{(0, i_0^{k_0})}\big\}_{k_0=1}^{n}$
\end{tabbing}

\begin{tabbing}
$\hat{c}_j := (1-\beta) c_j, \, j \in N$; \, 
$S_0^{0} := \emptyset$; \, $S_1^{0} := \emptyset$; \,
$k_0 := 0$ \\
\textbf{for} \= $k_1 := 1$ \textbf{to} $n$ \textbf{do} \\
\> $S_1^{k_1} := S_1^{k_1-1} \cup \{i_1^{k_1}\}$; \,
$\mathrm{AUGMENT}_1 := \texttt{false}$ \\
\> $\nu^{(0, k_1)}_{(0,j)} := \nu^{(k_1)}_j 
  - \hat{c}_j / w^{(k_1)}_j,  \, j \in S_1^{k_1} \setminus S_0^{k_0}$ \\
\> \textbf{while } \= $k_0 < k_1$ \textbf{ and } 
 \textbf{not}($\mathrm{AUGMENT}_1$) \textbf{ do } \\
\> \> \textbf{pick} 
 $j_0^{\max} \in \argmax \big\{\nu^{(0, k_1)}_{(0,j)}\colon j \in  S_1^{k_1} \setminus S_0^{k_0}\big\}$ \\ 
\> \> \textbf{if } \=  $k_1 = n$ \textbf{ or } $\nu^{*}_{i_1^{k_1}} < 
  \nu^{(0, k_1)}_{(0, j_0^{\max})}$ \\
\> \> \>  $i_0^{k_0+1} := j_0^{\max}$; \, $\nu^*_{(0, i_0^{k_0+1})} := 
\nu^{(0, k_1)}_{(0, i_0^{k_0+1})}$ \\
 \> \>  \>   $S_0^{k_0+1} := S_0^{k_0} \cup \{i_0^{k_0+1}\}$; \, $k_0 := k_0+1$  \\
\> \>  \textbf{else } \\
 \> \>  \>  $\mathrm{AUGMENT}_1 := \texttt{true}$ \\
 \> \> \textbf{end} \{ if \} \\
\> \textbf{end }  \{ while \} \\
\textbf{end }  \{ for \}
\end{tabbing}
\end{minipage}}
\end{center}
\label{fig:hatag3}
\end{table}

We next assess the computational complexity of the switching-index algorithm.

\begin{proposition}
\label{pro:cag0}
Algorithm $\mathrm{AG}^0$ performs at most $n^2 + O(n)$ arithmetic operations.
\end{proposition}
\begin{proof}
The operation count is dominated by the statements
\[
\nu^{(0, k_1)}_{(0,j)} := \nu^{(k_1)}_j 
  - \hat{c}_j / w^{(k_1)}_j,  \quad j \in S_1^{k_1} \setminus S_0^{k_0},
\]
performing at most $2 k_1$ operations. Adding up
over $k_1 = 1, \ldots, n$ yields the result.
\end{proof}

We thus see that the switching index can be efficiently computed \emph{an
order of magnitude faster} than the continuation index. 
Note also that algorithm $\mathrm{AG}^0$ involves handling
matrices of size $n \times n$, in contrast to the matrices of size  $2
n \times 2n$ required by the AT scheme, which further yields
substantial savings in expensive memory operations.

\section{Index Dependence on Switching Costs}
\label{s:sa}
We discuss next some properties on the indices' dependence on switching
costs, in the case $c_i \equiv c$ and $d_i \equiv d$.
We write below $\nu_{(1, i)}^*(d)$ --- as it does
not depend on  $c$ --- and 
$\nu_{(0, i)}^*(c, d)$, and 
denote by $\nu_i^*$ the Gittins index of the underlying
 project with no switching costs.

\begin{proposition}
\label{pro:swnui} \mbox{ }
\begin{itemize}
\item[\textup{(a)}] $\nu_{(1, i)}^*(d) = \nu_i^* + (1-\beta) d$.
\item[\textup{(b)}] For large enough $c + d$,
$\nu^*_{(0, i)}(c, d) = \nu^N_i - (1-\beta) c$.
\item[\textup{(c)}] $\nu_{(0, i)}^*(c, d)$ is piecewise linear  
convex in $(c, d)$, decreasing in $c$ and nonincreasing in $d$.
\end{itemize}
\end{proposition}
\begin{proof}
(a) Notice that
$\nu_{(1, i)}^*(d)$ is the Gittins index of a project with modified
rewards $\tilde{R}_j = R_j + (1-\beta) d$ (cf.\ 
Section \ref{s:rzscc}). The effect of such an addition of a constant term to 
rewards is to increment the Gittins index by the same
constant, which yields the result.

(b) Using the results in Sections \ref{s:rzscc} and  \ref{s:dpaidmpi} we obtain
\begin{equation}
\label{eq:nu0icd}
\nu_{(0, i)}^*(c, d) = \max_{S \subseteq N\colon i \in S} \frac{f_i^S - c -
  \big\{1 - (1-\beta) g_i^S\big\} d}{g_i^S} = 
(1-\beta) d + \max_{S \subseteq N\colon i \in S} \frac{f_i^S - (c+d)}{g_i^S},
\end{equation}
where $f_i^S$ is the reward measure of the underlying nonrestless project with 
rewards $R_j$ --- note that the corresponding reward
measure with modified rewards $\tilde{R}_j$ as above is $f_i^S(d)
= f_i^S + (1-\beta) d g_i^S$.
The second identity in (\ref{eq:nu0icd}) implies that, for $c + d$
large enough, term $(c+d)/g_i^S$ becomes dominant, and hence
the maximum value of the given expression is attained by
maximizing the denominator $g_i^S$. The latter's maximum value is achieved by $S=N$, for which
$g_i^N = 1/(1-\beta)$. Since $\nu_i^N = r_i^N/w_i^N = f_i^N/g_i^N$,
this yields the result.

(c) 
The first identity in (\ref{eq:nu0icd}) represents $\nu_{(0, i)}^*(c, d)$ as the
maximum of linear functions in $(c, d)$ that are decreasing in $c$ and
nonincreasing in $d$, which yields the result.
\end{proof}

Note that Proposition \ref{pro:swnui}(a)
shows that the incentive to stay on an active project
increases linearly with its shutdown cost. 
Also, as $\beta \nearrow 1$ the dependence on switching
costs vanishes since the switching index converges
to the undiscounted Gittins index.

\begin{figure}[ht!]
\centering
\begin{psfrags}%
\psfragscanon%
\psfrag{s01}[t][t]{\color[rgb]{0,0,0}\setlength{\tabcolsep}{0pt}\begin{tabular}{c}$c$\end{tabular}}%
\psfrag{s02}[b][b]{\color[rgb]{0,0,0}\setlength{\tabcolsep}{0pt}\begin{tabular}{c}$\nu_{(0, i)}^*(c, 0)$\end{tabular}}%
%
\psfrag{x01}[t][t]{0}%
\psfrag{x02}[t][t]{0.1}%
\psfrag{x03}[t][t]{0.2}%
\psfrag{x04}[t][t]{0.3}%
\psfrag{x05}[t][t]{0.4}%
\psfrag{x06}[t][t]{0.5}%
\psfrag{x07}[t][t]{0.6}%
\psfrag{x08}[t][t]{0.7}%
\psfrag{x09}[t][t]{0.8}%
\psfrag{x10}[t][t]{0.9}%
\psfrag{x11}[t][t]{1}%
\psfrag{x12}[t][r]{0}%
\psfrag{x13}[t][t]{}%
\psfrag{x14}[t][t]{}%
\psfrag{x15}[t][t]{}%
\psfrag{x16}[t][t]{}%
\psfrag{x17}[t][t]{}%
\psfrag{x18}[t][t]{}%
\psfrag{x19}[t][t]{}%
\psfrag{x20}[t][t]{}%
\psfrag{x21}[t][t]{}%
\psfrag{x22}[t][t]{2}%
%
\psfrag{v01}[r][r]{0}%
\psfrag{v02}[r][r]{0.1}%
\psfrag{v03}[r][r]{0.2}%
\psfrag{v04}[r][r]{0.3}%
\psfrag{v05}[r][r]{0.4}%
\psfrag{v06}[r][r]{0.5}%
\psfrag{v07}[r][r]{0.6}%
\psfrag{v08}[r][r]{0.7}%
\psfrag{v09}[r][r]{0.8}%
\psfrag{v10}[r][r]{0.9}%
\psfrag{v11}[r][r]{1}%
\psfrag{v12}[r][r]{0.4}%
\psfrag{v13}[r][r]{}%
\psfrag{v14}[r][r]{}%
\psfrag{v15}[r][r]{}%
\psfrag{v16}[r][r]{}%
\psfrag{v17}[r][r]{}%
\psfrag{v18}[r][r]{1}%
%
\includegraphics[height=2in,width=6.5in,keepaspectratio]{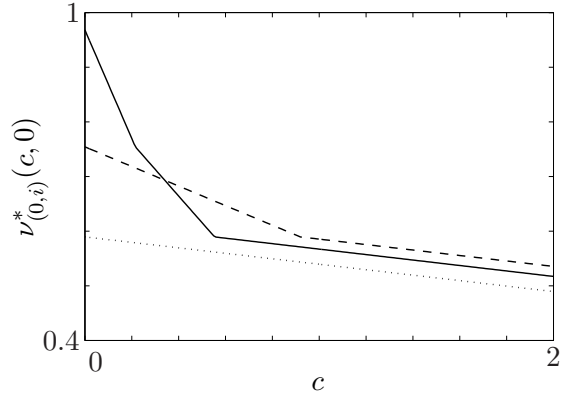}
\end{psfrags}%
\caption{Dependence of Switching Index on Startup Cost.}
\label{fig:mpi0suc}
\end{figure}

We next give two examples to illustrate the above results.
The first concerns the project instance with state space $N = \{1, 2,
3\}$, 
startup cost $c$,
\[
d = 0, \quad \beta = 0.95, \quad \mathbf{R}^0 =
\mathbf{0}, \quad \mathbf{R}^1 = \begin{bmatrix} 0.7221 \\ 0.9685 \\ 0.1557\end{bmatrix}, 
\quad \text{ and } \quad 
\mathbf{P} = 
\begin{bmatrix} 
0.8061 & 0.1574 & 0.0365 \\
0.1957 & 0.0067 & 0.7976 \\
0.1378 & 0.5959 & 0.2663
\end{bmatrix}.
\]
Figure \ref{fig:mpi0suc} plots  the 
switching index for each state vs.\  $c$, consistently
with Proposition
\ref{pro:swnui}(b, c), showing that the state
ordering induced by the index can change as $c$
varies.

\begin{figure}[ht!]
\centering
\begin{psfrags}%
\psfrag{s01}[t][t]{\color[rgb]{0,0,0}\setlength{\tabcolsep}{0pt}\begin{tabular}{c}$d$\end{tabular}}%
\psfrag{s02}[c][b]{\color[rgb]{0,0,0}\setlength{\tabcolsep}{0pt}\begin{tabular}{c}$\nu_{(1, i)}^*(d)$\end{tabular}}%
\psfrag{s05}[t][t]{\color[rgb]{0,0,0}\setlength{\tabcolsep}{0pt}\begin{tabular}{c}$d$\end{tabular}}%
\psfrag{s06}[c][b]{\color[rgb]{0,0,0}\setlength{\tabcolsep}{0pt}\begin{tabular}{c}$\nu_{(0,
      i)}^*(0, d)$\end{tabular}}%
%
\psfrag{x01}[t][t]{0}%
\psfrag{x02}[t][t]{0.1}%
\psfrag{x03}[t][t]{0.2}%
\psfrag{x04}[t][t]{0.3}%
\psfrag{x05}[t][t]{0.4}%
\psfrag{x06}[t][t]{0.5}%
\psfrag{x07}[t][t]{0.6}%
\psfrag{x08}[t][t]{0.7}%
\psfrag{x09}[t][t]{0.8}%
\psfrag{x10}[t][t]{0.9}%
\psfrag{x11}[t][t]{1}%
\psfrag{x12}[t][t]{0}%
\psfrag{x13}[t][t]{}%
\psfrag{x14}[t][t]{}%
\psfrag{x15}[t][t]{}%
\psfrag{x16}[t][t]{}%
\psfrag{x17}[t][t]{}%
\psfrag{x18}[t][t]{}%
\psfrag{x19}[t][t]{}%
\psfrag{x20}[t][t]{}%
\psfrag{x21}[t][t]{}%
\psfrag{x22}[t][t]{2}%
\psfrag{x23}[t][t]{0}%
\psfrag{x24}[t][t]{}%
\psfrag{x25}[t][t]{}%
\psfrag{x26}[t][t]{}%
\psfrag{x27}[t][t]{}%
\psfrag{x28}[t][t]{}%
\psfrag{x29}[t][t]{}%
\psfrag{x30}[t][t]{}%
\psfrag{x31}[t][t]{}%
\psfrag{x32}[t][t]{}%
\psfrag{x33}[t][t]{2}%
%
\psfrag{v01}[r][r]{0}%
\psfrag{v02}[r][r]{0.1}%
\psfrag{v03}[r][r]{0.2}%
\psfrag{v04}[r][r]{0.3}%
\psfrag{v05}[r][r]{0.4}%
\psfrag{v06}[r][r]{0.5}%
\psfrag{v07}[r][r]{0.6}%
\psfrag{v08}[r][r]{0.7}%
\psfrag{v09}[r][r]{0.8}%
\psfrag{v10}[r][r]{0.9}%
\psfrag{v11}[r][r]{1}%
\psfrag{v12}[r][r]{0.5}%
\psfrag{v13}[r][r]{}%
\psfrag{v14}[r][r]{}%
\psfrag{v15}[r][r]{}%
\psfrag{v16}[r][r]{}%
\psfrag{v17}[r][r]{}%
\psfrag{v18}[r][r]{1.1}%
\psfrag{v19}[r][r]{0.5}%
\psfrag{v20}[r][r]{}%
\psfrag{v21}[r][r]{}%
\psfrag{v22}[r][r]{}%
\psfrag{v23}[r][r]{}%
\psfrag{v24}[r][r]{}%
\psfrag{v25}[r][r]{1.1}%
\includegraphics[height=2in,width=6.5in]{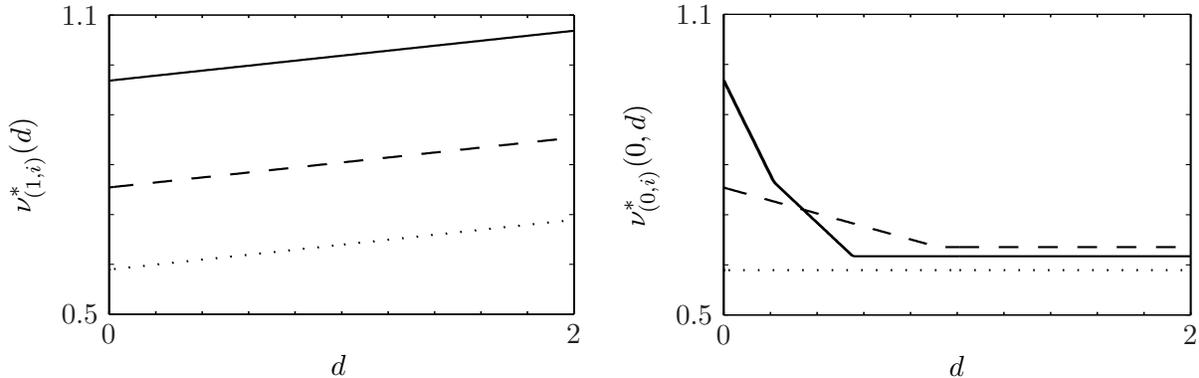}
\end{psfrags}%
\caption{Dependence of Continuation and Switching Indices on Shutdown Cost.}
\label{fig:mpisdc}
\end{figure}

In the second example we let the
shutdown cost $d$ vary, taking $c = 0$.
Figure \ref{fig:mpisdc} plots  the
indices for each state vs.\ $d$.
Again, the plots are consistent with Proposition
\ref{pro:swnui}. 

\section{Computational Experiments}
\label{s:es}
We next report on a computational study 
based on the author's implementations of the results
herein.
The experiments
were performed running MATLAB R2006b under Windows XP x64 in an
HP xw9300 AMD Opteron 254 (2.8 GHz) workstation with 4 GB of memory.

The first experiment
investigated the runtime performance of the 
decoupled index computation method. 
We made MATLAB generate
 a random project instance with startup costs for each 
of the state-space sizes $n = 500, 1000, \ldots, 5000$.
For each such $n$, MATLAB recorded the 
time  to compute: (i) the continuation index and 
additional quantities with algorithm FP(1)
 in \citet{nmijoc06}; (ii) the switching index by
algorithm $\mathrm{AG}^0$;  and (iii) both
indices via the AT scheme --- using algorithm FP(0) in \citet{nmijoc06}.

\begin{figure}[ht!]
\centering
\begin{psfrags}%
\psfragscanon%
\psfrag{s01}[t][t]{\color[rgb]{0,0,0}\setlength{\tabcolsep}{0pt}\begin{tabular}{c}$n$\end{tabular}}%
\psfrag{s02}[c][b]{\color[rgb]{0,0,0}\setlength{\tabcolsep}{0pt}\begin{tabular}{c}runtime (hours)\end{tabular}}%
\psfrag{s03}[t][t]{\color[rgb]{0,0,0}\setlength{\tabcolsep}{0pt}\begin{tabular}{c}AT scheme\end{tabular}}%
\psfrag{s04}[t][t]{\color[rgb]{0,0,0}\setlength{\tabcolsep}{0pt}\begin{tabular}{c}two-stage method\end{tabular}}%
\psfrag{s07}[t][t]{\color[rgb]{0,0,0}\setlength{\tabcolsep}{0pt}\begin{tabular}{c}$n$\end{tabular}}%
\psfrag{s08}[c][b]{\color[rgb]{0,0,0}\setlength{\tabcolsep}{0pt}\begin{tabular}{c}runtime (secs.) of $\mathrm{AG}^0$\end{tabular}}%
%
\psfrag{x01}[t][t]{0}%
\psfrag{x02}[t][t]{0.1}%
\psfrag{x03}[t][t]{0.2}%
\psfrag{x04}[t][t]{0.3}%
\psfrag{x05}[t][t]{0.4}%
\psfrag{x06}[t][t]{0.5}%
\psfrag{x07}[t][t]{0.6}%
\psfrag{x08}[t][t]{0.7}%
\psfrag{x09}[t][t]{0.8}%
\psfrag{x10}[t][t]{0.9}%
\psfrag{x11}[t][t]{1}%
\psfrag{x12}[t][t]{500}%
\psfrag{x13}[t][t]{}%
\psfrag{x14}[t][t]{}%
\psfrag{x15}[t][t]{}%
\psfrag{x16}[t][t]{}%
\psfrag{x17}[t][t]{}%
\psfrag{x18}[t][t]{}%
\psfrag{x19}[t][t]{}%
\psfrag{x20}[t][t]{}%
\psfrag{x21}[t][t]{5000}%
\psfrag{x22}[t][t]{500}%
\psfrag{x23}[t][t]{}%
\psfrag{x24}[t][t]{}%
\psfrag{x25}[t][t]{}%
\psfrag{x26}[t][t]{}%
\psfrag{x27}[t][t]{}%
\psfrag{x28}[t][t]{}%
\psfrag{x29}[t][t]{}%
\psfrag{x30}[t][t]{}%
\psfrag{x31}[t][t]{5000}%
%
\psfrag{v01}[r][r]{0}%
\psfrag{v02}[r][r]{0.1}%
\psfrag{v03}[r][r]{0.2}%
\psfrag{v04}[r][r]{0.3}%
\psfrag{v05}[r][r]{0.4}%
\psfrag{v06}[r][r]{0.5}%
\psfrag{v07}[r][r]{0.6}%
\psfrag{v08}[r][r]{0.7}%
\psfrag{v09}[r][r]{0.8}%
\psfrag{v10}[r][r]{0.9}%
\psfrag{v11}[r][r]{1}%
\psfrag{v12}[r][r]{0}%
\psfrag{v13}[r][r]{}%
\psfrag{v14}[r][r]{1}%
\psfrag{v15}[r][r]{}%
\psfrag{v16}[r][r]{2}%
\psfrag{v17}[r][r]{}%
\psfrag{v18}[r][r]{3}%
\psfrag{v19}[r][r]{}%
\psfrag{v20}[r][r]{4}%
\psfrag{v21}[r][r]{0}%
\psfrag{v22}[r][r]{}%
\psfrag{v23}[r][r]{1}%
\psfrag{v24}[r][r]{}%
\psfrag{v25}[r][r]{2}%
\psfrag{v26}[r][r]{}%
\psfrag{v27}[r][r]{3}%
\psfrag{v28}[r][r]{}%
\psfrag{v29}[r][r]{4}%
%
\includegraphics[height=2.8in,width=6.5in,keepaspectratio]{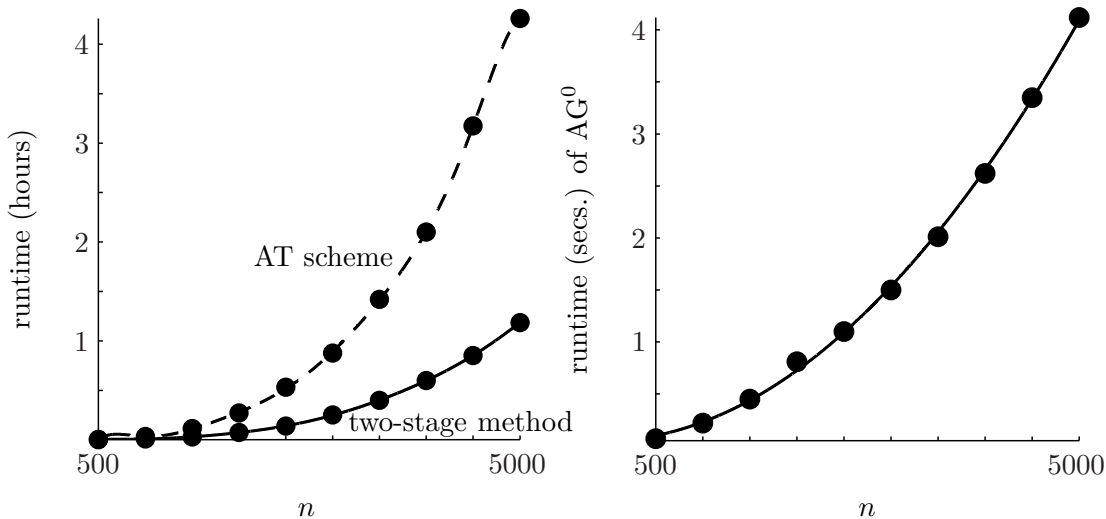}
\end{psfrags}%
\caption{Exp.\ 1: Runtimes of Index Algorithms.}
\label{fig:plot1_timingexp}
\end{figure}

The results are displayed in Figure \ref{fig:plot1_timingexp}. 
The left pane plots total runtimes in hours for computing both
indices vs.\ $n$, along with cubic least-squares (LS)
fits, which are consistent with the theoretical $O(n^3)$ complexity.
The dotted line corresponds to the AT scheme, while the solid line corresponds to our two-stage method. 
The results show that the latter consistently achieved about a 
four-fold speedup over the former.

The right pane plots runtimes, in \emph{seconds}, for the switching-index
algorithm vs.\ $n$, along with a quadratic least-squares
fit, which is consistent with the theoretical $O(n^2)$ complexity.
The change of timescale demonstrates the
order-of-magnitude runtime improvement.

In the following experiments 
 the optimal policy was 
computed for each instance solving the linear programming (LP) formulation of the DP equations 
using the CPLEX LP solver, interfaced with MATLAB via TOMLAB.
The MPI and
benchmark (Gittins index) policies were evaluated solving with MATLAB their evaluation equations.

The second experiment assessed how the
performance of the MPI policy on two-project
instances  depends on a common constant startup cost $c$ and
discount factor $\beta$ --- shutdown costs are zero.
A sample of
100 instances (where each project has state space $N = \{1, \ldots, n\}$ with $n = 10$) 
was randomly generated.
In each instance, parameter values for each project were independently
generated: transition probabilities (obtained by scaling a matrix with
Uniform[0, 1] entries  dividing each row by its sum) and active
rewards (Uniform[0, 1]). Passive rewards were set to zero.
For each instance $k = 1, \ldots, 100$ and pair
$(c, \beta) \in [0, 1] \times [0.2, 0.9]$ --- using a 0.1 grid ---
the optimal objective value $v^{(k), \textup{opt}}$ and the
objective values  of the MPI ($v^{(k), \textup{MPI}}$) and the benchmark ($v^{(k), \textup{bench}}$)
policies were computed, along with the corresponding 
relative suboptimality gap of the MPI policy
$\Delta^{(k), \textup{MPI}}
 = 100 (v^{(k), \textup{opt}} - v^{(k), \textup{MPI}}) / |v^{(k),
    \textup{opt}}|$,
and the suboptimality-gap ratio of the MPI over the 
benchmark policy $\rho^{(k), \textup{MPI}, \textup{bench}} = 100
(v^{(k), \textup{opt}} - v^{(k), \textup{MPI}}) / 
(v^{(k), \textup{opt}} - v^{(k), \textup{bench}})$.
The latter were 
then averaged over the 100 instances for each $(c, \beta)$ pair to
obtain the average values $\Delta^{\textup{MPI}}$ and $\rho^{\textup{MPI}, \textup{bench}}$.

Objective values $v^{(k), \textup{opt}}$, $v^{(k),
  \textup{MPI}}$ and $v^{(k), \textup{bench}}$ were evaluated as follows.
First, the corresponding \emph{value functions} $v^{(k), \textup{opt}}_{((a_1^-,
  i_1), (a_2^-, i_2))}$, $v^{(k),
  \textup{MPI}}_{((a_1^-,
  i_1), (a_2^-, i_2))}$ and $v^{(k), \textup{bench}}_{((a_1^-,
  i_1), (a_2^-, i_2))}$
were computed as mentioned above, as vectors indexed by the initial
  joint augmented state $((a_1^-,
  i_1), (a_2^-, i_2))$.
Then, the objective values were evaluated by taking both projects to be initially passive:
\begin{equation}
\label{eq:vkpi}
v^{(k), \pi} \triangleq \frac{1}{n^2} \sum_{i_1, i_2 \in N} v^{(k),
  \pi}_{((0, i_1), (0, i_2))}, \quad \pi \in \{\textup{opt},
  \textup{MPI}, \textup{bench}\}.
\end{equation}

\begin{figure}[ht!]
\centering
\begin{psfrags}%
\psfragscanon%
\psfrag{s01}[c][t]{\color[rgb]{0,0,0}\setlength{\tabcolsep}{0pt}\begin{tabular}{c}$c$\end{tabular}}%
\psfrag{s02}[t][b]{\color[rgb]{0,0,0}\setlength{\tabcolsep}{0pt}\begin{tabular}{c}$\Delta^{\textup{MPI}}$\end{tabular}}%
\psfrag{s03}[b][b]{\color[rgb]{0,0,0}\setlength{\tabcolsep}{0pt}\begin{tabular}{c}Dependence on $c$ for Multiple $\beta$\end{tabular}}%
\psfrag{s05}[c][t]{\color[rgb]{0,0,0}\setlength{\tabcolsep}{0pt}\begin{tabular}{c}$\beta$\end{tabular}}%
\psfrag{s06}[t][b]{\color[rgb]{0,0,0}\setlength{\tabcolsep}{0pt}\begin{tabular}{c}$\Delta^{\textup{MPI}}$\end{tabular}}%
\psfrag{s07}[b][b]{\color[rgb]{0,0,0}\setlength{\tabcolsep}{0pt}\begin{tabular}{c}Dependence on $\beta$ for $c = 0.3$\end{tabular}}%
%
\psfrag{x01}[t][t]{0}%
\psfrag{x02}[t][t]{0.1}%
\psfrag{x03}[t][t]{0.2}%
\psfrag{x04}[t][t]{0.3}%
\psfrag{x05}[t][t]{0.4}%
\psfrag{x06}[t][t]{0.5}%
\psfrag{x07}[t][t]{0.6}%
\psfrag{x08}[t][t]{0.7}%
\psfrag{x09}[t][t]{0.8}%
\psfrag{x10}[t][t]{0.9}%
\psfrag{x11}[t][t]{1}%
\psfrag{x12}[t][r]{0.2}%
\psfrag{x13}[t][t]{}%
\psfrag{x14}[t][t]{}%
\psfrag{x15}[t][t]{}%
\psfrag{x16}[t][t]{}%
\psfrag{x17}[t][t]{}%
\psfrag{x18}[t][t]{}%
\psfrag{x19}[t][t]{0.9}%
\psfrag{x20}[t][t]{0}%
\psfrag{x21}[t][t]{}%
\psfrag{x22}[t][t]{}%
\psfrag{x23}[t][t]{}%
\psfrag{x24}[t][t]{}%
\psfrag{x25}[t][t]{}%
\psfrag{x26}[t][t]{}%
\psfrag{x27}[t][t]{}%
\psfrag{x28}[t][t]{}%
\psfrag{x29}[t][t]{}%
\psfrag{x30}[t][t]{1}%
%
\psfrag{v01}[r][r]{0}%
\psfrag{v02}[r][r]{0.1}%
\psfrag{v03}[r][r]{0.2}%
\psfrag{v04}[r][r]{0.3}%
\psfrag{v05}[r][r]{0.4}%
\psfrag{v06}[r][r]{0.5}%
\psfrag{v07}[r][r]{0.6}%
\psfrag{v08}[r][r]{0.7}%
\psfrag{v09}[r][r]{0.8}%
\psfrag{v10}[r][r]{0.9}%
\psfrag{v11}[r][r]{1}%
\psfrag{v12}[r][r]{0\%}%
\psfrag{v13}[r][r]{}%
\psfrag{v14}[r][r]{}%
\psfrag{v15}[r][r]{}%
\psfrag{v16}[r][r]{}%
\psfrag{v17}[r][r]{0.25\%}%
\psfrag{v18}[r][r]{0\%}%
\psfrag{v19}[r][r]{}%
\psfrag{v20}[r][r]{}%
\psfrag{v21}[r][r]{}%
\psfrag{v22}[r][r]{}%
\psfrag{v23}[r][r]{0.25\%}%
%
\includegraphics[height=2in,width=6.5in]{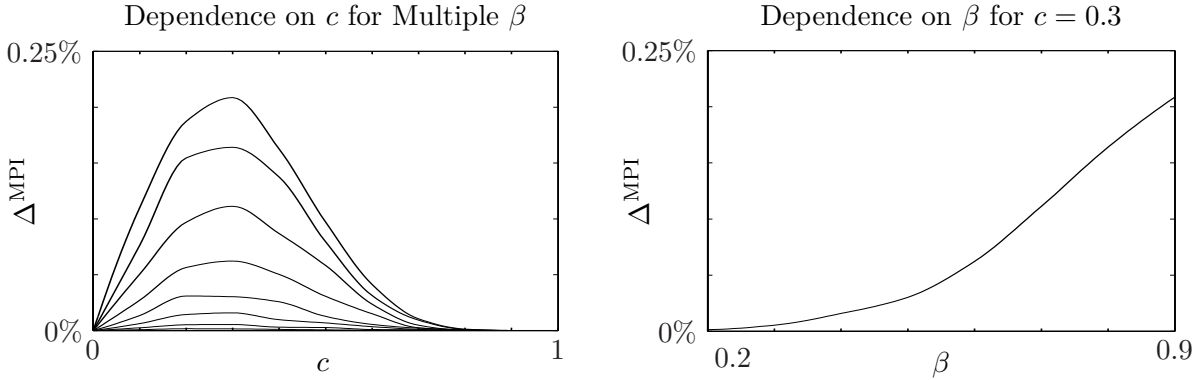}
\end{psfrags}%
\caption{Exp.\ 2: Average Relative Suboptimality Gap of MPI Policy.}
\label{fig:exp1aplot}
\end{figure}

The left pane in Figure \ref{fig:exp1aplot} plots $\Delta^{\textup{MPI}}$ vs.\ the
 startup cost $c$ for multiple $\beta$, using MATLAB's
 cubic interpolation.
Such a gap starts at $0$ for $c = 0$
--- as  the optimal policy is then
recovered ---
then increases up to a maximum value below $0.25\%$ at around $c \approx 0.3$, 
and then decreases, dropping again to $0$ at about $c \approx 0.9$
and staying there for larger values of $c$.
Such a behavior is to 
be expected: for large enough $c$ the optimal
 policy will simply pick a project and hold to it. 
The plot in the right pane shows that 
$\Delta^{\textup{MPI}}$ increases with $\beta$.
 
\begin{figure}[ht!]
\centering
\begin{psfrags}%
\psfragscanon%
\psfrag{s01}[c][t]{\color[rgb]{0,0,0}\setlength{\tabcolsep}{0pt}\begin{tabular}{c}$c$\end{tabular}}%
\psfrag{s02}[t][b]{\color[rgb]{0,0,0}\setlength{\tabcolsep}{0pt}\begin{tabular}{c}$\rho^{\textup{MPI}, \textup{bench}}$\end{tabular}}%
\psfrag{s03}[b][b]{\color[rgb]{0,0,0}\setlength{\tabcolsep}{0pt}\begin{tabular}{c}Dependence on $c$ for Multiple $\beta$\end{tabular}}%
\psfrag{s05}[t][t]{\color[rgb]{0,0,0}\setlength{\tabcolsep}{0pt}\begin{tabular}{c}$\beta$\end{tabular}}%
\psfrag{s06}[t][b]{\color[rgb]{0,0,0}\setlength{\tabcolsep}{0pt}\begin{tabular}{c}$\rho^{\textup{MPI}, \textup{bench}}$\end{tabular}}%
\psfrag{s07}[b][b]{\color[rgb]{0,0,0}\setlength{\tabcolsep}{0pt}\begin{tabular}{c}Dependence on $\beta$ for $c = 0.1$\end{tabular}}%
%
\psfrag{x01}[t][t]{0}%
\psfrag{x02}[t][t]{0.1}%
\psfrag{x03}[t][t]{0.2}%
\psfrag{x04}[t][t]{0.3}%
\psfrag{x05}[t][t]{0.4}%
\psfrag{x06}[t][t]{0.5}%
\psfrag{x07}[t][t]{0.6}%
\psfrag{x08}[t][t]{0.7}%
\psfrag{x09}[t][t]{0.8}%
\psfrag{x10}[t][t]{0.9}%
\psfrag{x11}[t][t]{1}%
\psfrag{x12}[t][r]{0.2}%
\psfrag{x13}[t][t]{}%
\psfrag{x14}[t][t]{}%
\psfrag{x15}[t][t]{}%
\psfrag{x16}[t][t]{}%
\psfrag{x17}[t][t]{}%
\psfrag{x18}[t][t]{}%
\psfrag{x19}[t][t]{0.9}%
\psfrag{x20}[t][t]{0}%
\psfrag{x21}[t][t]{}%
\psfrag{x22}[t][t]{}%
\psfrag{x23}[t][t]{}%
\psfrag{x24}[t][t]{}%
\psfrag{x25}[t][t]{}%
\psfrag{x26}[t][t]{}%
\psfrag{x27}[t][t]{}%
\psfrag{x28}[t][t]{}%
\psfrag{x29}[t][t]{}%
\psfrag{x30}[t][t]{1}%
%
\psfrag{v01}[r][r]{0}%
\psfrag{v02}[r][r]{0.1}%
\psfrag{v03}[r][r]{0.2}%
\psfrag{v04}[r][r]{0.3}%
\psfrag{v05}[r][r]{0.4}%
\psfrag{v06}[r][r]{0.5}%
\psfrag{v07}[r][r]{0.6}%
\psfrag{v08}[r][r]{0.7}%
\psfrag{v09}[r][r]{0.8}%
\psfrag{v10}[r][r]{0.9}%
\psfrag{v11}[r][r]{1}%
\psfrag{v12}[r][r]{0\%}%
\psfrag{v13}[r][r]{}%
\psfrag{v14}[r][r]{}%
\psfrag{v15}[r][r]{}%
\psfrag{v16}[r][r]{}%
\psfrag{v17}[r][r]{}%
\psfrag{v18}[r][r]{}%
\psfrag{v19}[r][r]{}%
\psfrag{v20}[r][r]{40\%}%
\psfrag{v21}[r][r]{0\%}%
\psfrag{v22}[r][r]{}%
\psfrag{v23}[r][r]{}%
\psfrag{v24}[r][r]{}%
\psfrag{v25}[r][r]{}%
\psfrag{v26}[r][r]{}%
\psfrag{v27}[r][r]{}%
\psfrag{v28}[r][r]{}%
\psfrag{v29}[r][r]{40\%}%
%
\includegraphics[height=2in,width=6.5in]{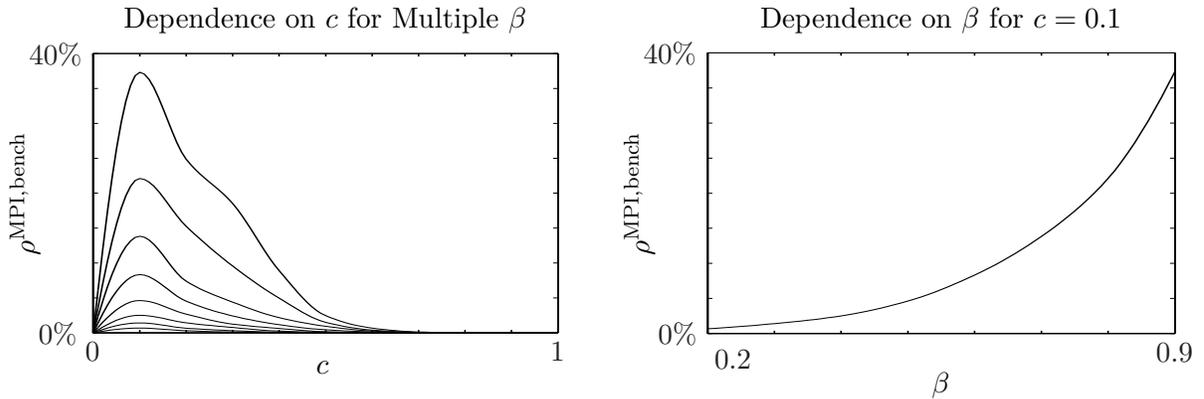}
\end{psfrags}%
\caption{Exp.\ 2: Average Suboptimality-Gap Ratio of MPI over Benchmark
  Policy.}
\label{fig:exp1bplot}
\end{figure}

Figure \ref{fig:exp1bplot} displays corresponding plots for the suboptimality-gap
ratio $\rho^{\textup{MPI}, \textup{bench}}$ of the MPI over the benchmark policy. 
The plot in the left pane shows that the average suboptimality gap for
the MPI policy remains below $40\%$ of the gap for the benchmark
policy. Further, such a ratio takes the value $0$ both for $c = 0$ and for
$c$ large enough, as the MPI policy is then optimal. The plot in the
right pane shows that the
ratio increases with $\beta$. 

\begin{figure}[ht!]
\centering
\begin{psfrags}%
\psfragscanon%
\psfrag{s01}[c][t]{\color[rgb]{0,0,0}\setlength{\tabcolsep}{0pt}\begin{tabular}{c}$d$\end{tabular}}%
\psfrag{s02}[t][b]{\color[rgb]{0,0,0}\setlength{\tabcolsep}{0pt}\begin{tabular}{c}$\Delta^{\textup{MPI}}$\end{tabular}}%
\psfrag{s03}[b][b]{\color[rgb]{0,0,0}\setlength{\tabcolsep}{0pt}\begin{tabular}{c}Dependence on $d$ for Multiple $\beta$\end{tabular}}%
\psfrag{s05}[c][t]{\color[rgb]{0,0,0}\setlength{\tabcolsep}{0pt}\begin{tabular}{c}$\beta$\end{tabular}}%
\psfrag{s06}[t][b]{\color[rgb]{0,0,0}\setlength{\tabcolsep}{0pt}\begin{tabular}{c}$\Delta^{\textup{MPI}}$\end{tabular}}%
\psfrag{s07}[b][b]{\color[rgb]{0,0,0}\setlength{\tabcolsep}{0pt}\begin{tabular}{c}Dependence on $\beta$ for $d = 0.3$\end{tabular}}%
%
\psfrag{x01}[t][t]{0}%
\psfrag{x02}[t][t]{0.1}%
\psfrag{x03}[t][r]{0.2}%
\psfrag{x04}[t][t]{0.3}%
\psfrag{x05}[t][t]{0.4}%
\psfrag{x06}[t][t]{0.5}%
\psfrag{x07}[t][t]{0.6}%
\psfrag{x08}[t][t]{0.7}%
\psfrag{x09}[t][t]{0.8}%
\psfrag{x10}[t][t]{0.9}%
\psfrag{x11}[t][t]{1}%
\psfrag{x12}[t][r]{0.2}%
\psfrag{x13}[t][t]{}%
\psfrag{x14}[t][t]{}%
\psfrag{x15}[t][t]{}%
\psfrag{x16}[t][t]{}%
\psfrag{x17}[t][t]{}%
\psfrag{x18}[t][t]{}%
\psfrag{x19}[t][t]{0.9}%
\psfrag{x20}[t][t]{0}%
\psfrag{x21}[t][t]{}%
\psfrag{x22}[t][t]{}%
\psfrag{x23}[t][t]{}%
\psfrag{x24}[t][t]{}%
\psfrag{x25}[t][t]{}%
\psfrag{x26}[t][t]{}%
\psfrag{x27}[t][t]{}%
\psfrag{x28}[t][t]{}%
\psfrag{x29}[t][t]{}%
\psfrag{x30}[t][t]{1}%
%
\psfrag{v01}[r][r]{0}%
\psfrag{v02}[r][r]{0.1}%
\psfrag{v03}[r][r]{0.2}%
\psfrag{v04}[r][r]{0.3}%
\psfrag{v05}[r][r]{0.4}%
\psfrag{v06}[r][r]{0.5}%
\psfrag{v07}[r][r]{0.6}%
\psfrag{v08}[r][r]{0.7}%
\psfrag{v09}[r][r]{0.8}%
\psfrag{v10}[r][r]{0.9}%
\psfrag{v11}[r][r]{1}%
\psfrag{v12}[r][r]{0\%}%
\psfrag{v13}[r][r]{}%
\psfrag{v14}[r][r]{}%
\psfrag{v15}[r][r]{}%
\psfrag{v16}[r][r]{}%
\psfrag{v17}[r][r]{0.25\%}%
\psfrag{v18}[r][r]{0\%}%
\psfrag{v19}[r][r]{}%
\psfrag{v20}[r][r]{}%
\psfrag{v21}[r][r]{}%
\psfrag{v22}[r][r]{}%
\psfrag{v23}[r][r]{0.25\%}%
\includegraphics[height=2in,width=6.5in]{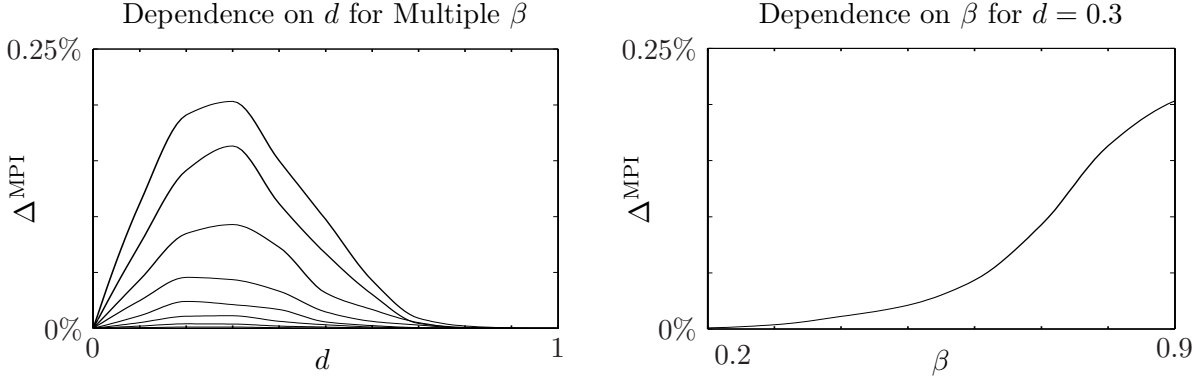}
\end{psfrags}%
\caption{Exp.\ 3: Average Relative Suboptimality Gap of MPI Policy.}
\label{fig:exp2aplot}
\end{figure}

\begin{figure}[ht!]
\centering
\begin{psfrags}%
\psfragscanon%
\psfrag{s01}[c][t]{\color[rgb]{0,0,0}\setlength{\tabcolsep}{0pt}\begin{tabular}{c}$d$\end{tabular}}%
\psfrag{s02}[t][b]{\color[rgb]{0,0,0}\setlength{\tabcolsep}{0pt}\begin{tabular}{c}$\rho^{\textup{MPI}, \textup{bench}}$\end{tabular}}%
\psfrag{s03}[b][b]{\color[rgb]{0,0,0}\setlength{\tabcolsep}{0pt}\begin{tabular}{c}Dependence on $d$ for Multiple $\beta$\end{tabular}}%
\psfrag{s05}[c][t]{\color[rgb]{0,0,0}\setlength{\tabcolsep}{0pt}\begin{tabular}{c}$\beta$\end{tabular}}%
\psfrag{s06}[t][b]{\color[rgb]{0,0,0}\setlength{\tabcolsep}{0pt}\begin{tabular}{c}$\rho^{\textup{MPI}, \textup{bench}}$\end{tabular}}%
\psfrag{s07}[b][b]{\color[rgb]{0,0,0}\setlength{\tabcolsep}{0pt}\begin{tabular}{c}Dependence on $\beta$ for $d = 0.1$\end{tabular}}%
%
\psfrag{x01}[t][t]{0}%
\psfrag{x02}[t][t]{0.1}%
\psfrag{x03}[t][t]{0.2}%
\psfrag{x04}[t][t]{0.3}%
\psfrag{x05}[t][t]{0.4}%
\psfrag{x06}[t][t]{0.5}%
\psfrag{x07}[t][t]{0.6}%
\psfrag{x08}[t][t]{0.7}%
\psfrag{x09}[t][t]{0.8}%
\psfrag{x10}[t][t]{0.9}%
\psfrag{x11}[t][t]{1}%
\psfrag{x12}[t][r]{0.2}%
\psfrag{x13}[t][t]{}%
\psfrag{x14}[t][t]{}%
\psfrag{x15}[t][t]{}%
\psfrag{x16}[t][t]{}%
\psfrag{x17}[t][t]{}%
\psfrag{x18}[t][t]{}%
\psfrag{x19}[t][t]{0.9}%
\psfrag{x20}[t][t]{0}%
\psfrag{x21}[t][t]{}%
\psfrag{x22}[t][t]{}%
\psfrag{x23}[t][t]{}%
\psfrag{x24}[t][t]{}%
\psfrag{x25}[t][t]{}%
\psfrag{x26}[t][t]{}%
\psfrag{x27}[t][t]{}%
\psfrag{x28}[t][t]{}%
\psfrag{x29}[t][t]{}%
\psfrag{x30}[t][t]{1}%
%
\psfrag{v01}[r][r]{0}%
\psfrag{v02}[r][r]{0.1}%
\psfrag{v03}[r][r]{0.2}%
\psfrag{v04}[r][r]{0.3}%
\psfrag{v05}[r][r]{0.4}%
\psfrag{v06}[r][r]{0.5}%
\psfrag{v07}[r][r]{0.6}%
\psfrag{v08}[r][r]{0.7}%
\psfrag{v09}[r][r]{0.8}%
\psfrag{v10}[r][r]{0.9}%
\psfrag{v11}[r][r]{1}%
\psfrag{v12}[r][r]{0\%}%
\psfrag{v13}[r][r]{}%
\psfrag{v14}[r][r]{}%
\psfrag{v15}[r][r]{}%
\psfrag{v16}[r][r]{}%
\psfrag{v17}[r][r]{}%
\psfrag{v18}[r][r]{}%
\psfrag{v19}[r][r]{35\%}%
\psfrag{v20}[r][r]{0\%}%
\psfrag{v21}[r][r]{}%
\psfrag{v22}[r][r]{}%
\psfrag{v23}[r][r]{}%
\psfrag{v24}[r][r]{}%
\psfrag{v25}[r][r]{}%
\psfrag{v26}[r][r]{}%
\psfrag{v27}[r][r]{35\%}%
%
\includegraphics[height=2in,width=6.5in]{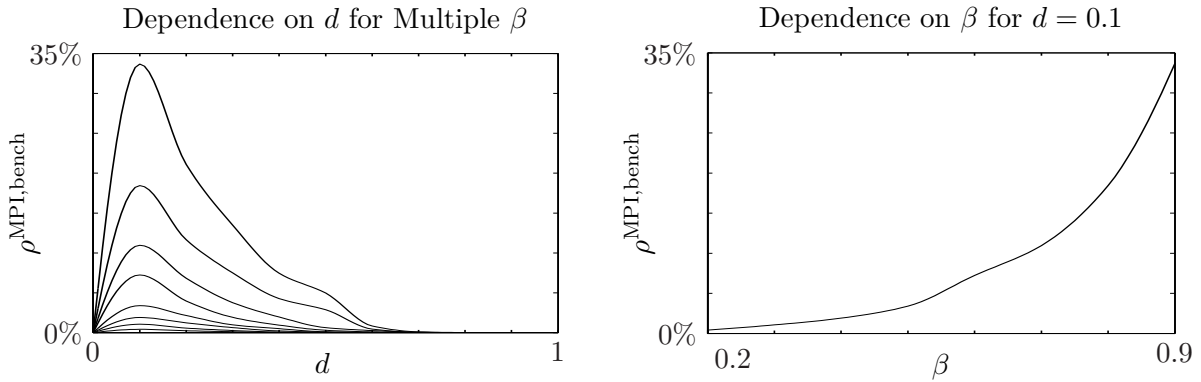}
\end{psfrags}%
\caption{Exp.\ 3: Average Suboptimality-Gap Ratio of MPI over Benchmark
  Policy.}
\label{fig:exp2bplot}
\end{figure}

The third experiment investigated the influence of a common
shutdown cost $d$ on the relative performance of the MPI policy ---
startup costs were set to zero. It was identical in other respects to the
second experiment.
The resultant plots, shown in
Figures \ref{fig:exp2aplot} and \ref{fig:exp2bplot}, are similar 
to those of the previous experiment.
Hence, the influence of a common shutdown cost on relative performance 
appears to be similar to that of a common startup cost. 

\begin{figure}[ht!]
\centering
\begin{psfrags}%
\psfragscanon%
\psfrag{s01}[t][t]{\color[rgb]{0,0,0}\setlength{\tabcolsep}{0pt}\begin{tabular}{c}$c_1$\end{tabular}}%
\psfrag{s02}[t][b]{\color[rgb]{0,0,0}\setlength{\tabcolsep}{0pt}\begin{tabular}{c}$c_2$\end{tabular}}%
\psfrag{s04}[b][b]{\color[rgb]{0,0,0}\setlength{\tabcolsep}{0pt}\begin{tabular}{c}$\Delta^{\textup{MPI}}$\end{tabular}}%
\psfrag{s05}[][]{\color[rgb]{0,0,0}\setlength{\tabcolsep}{0pt}\begin{tabular}{c}$0.04$\end{tabular}}%
\psfrag{s06}[][]{\color[rgb]{0,0,0}\setlength{\tabcolsep}{0pt}\begin{tabular}{c}$0.08$\end{tabular}}%
\psfrag{s07}[][]{\color[rgb]{0,0,0}\setlength{\tabcolsep}{0pt}\begin{tabular}{c}$0.12$\end{tabular}}%
\psfrag{s08}[][]{\color[rgb]{0,0,0}\setlength{\tabcolsep}{0pt}\begin{tabular}{c}$0.16$\end{tabular}}%
\psfrag{s09}[][]{\color[rgb]{0,0,0}\setlength{\tabcolsep}{0pt}\begin{tabular}{c}$0.2$\end{tabular}}%
\psfrag{s10}[t][t]{\color[rgb]{0,0,0}\setlength{\tabcolsep}{0pt}\begin{tabular}{c}$c_1$\end{tabular}}%
\psfrag{s11}[t][b]{\color[rgb]{0,0,0}\setlength{\tabcolsep}{0pt}\begin{tabular}{c}$c_2$\end{tabular}}%
\psfrag{s13}[b][b]{\color[rgb]{0,0,0}\setlength{\tabcolsep}{0pt}\begin{tabular}{c}$\rho^{\textup{MPI}, \textup{bench}}$\end{tabular}}%
\psfrag{s14}[][]{\color[rgb]{0,0,0}\setlength{\tabcolsep}{0pt}\begin{tabular}{c}$35$\end{tabular}}%
\psfrag{s15}[][]{\color[rgb]{0,0,0}\setlength{\tabcolsep}{0pt}\begin{tabular}{c}$28$\end{tabular}}%
\psfrag{s16}[][]{\color[rgb]{0,0,0}\setlength{\tabcolsep}{0pt}\begin{tabular}{c}$21$\end{tabular}}%
\psfrag{s17}[][]{\color[rgb]{0,0,0}\setlength{\tabcolsep}{0pt}\begin{tabular}{c}$14$\end{tabular}}%
\psfrag{s18}[][]{\color[rgb]{0,0,0}\setlength{\tabcolsep}{0pt}\begin{tabular}{c}$7$\end{tabular}}%
\psfrag{s19}[][]{\color[rgb]{0,0,0}\setlength{\tabcolsep}{0pt}\begin{tabular}{c}$1$\end{tabular}}%
%
\psfrag{x01}[t][t]{0}%
\psfrag{x02}[t][t]{0.1}%
\psfrag{x03}[t][t]{0.2}%
\psfrag{x04}[t][t]{0.3}%
\psfrag{x05}[t][t]{0.4}%
\psfrag{x06}[t][t]{0.5}%
\psfrag{x07}[t][t]{0.6}%
\psfrag{x08}[t][t]{0.7}%
\psfrag{x09}[t][t]{0.8}%
\psfrag{x10}[t][t]{0.9}%
\psfrag{x11}[t][t]{1}%
\psfrag{x12}[t][t]{0}%
\psfrag{x13}[t][t]{}%
\psfrag{x14}[t][t]{}%
\psfrag{x15}[t][t]{}%
\psfrag{x16}[t][t]{}%
\psfrag{x17}[t][t]{}%
\psfrag{x18}[t][t]{}%
\psfrag{x19}[t][t]{}%
\psfrag{x20}[t][t]{0.8}%
\psfrag{x21}[t][t]{0}%
\psfrag{x22}[t][t]{}%
\psfrag{x23}[t][t]{}%
\psfrag{x24}[t][t]{}%
\psfrag{x25}[t][t]{}%
\psfrag{x26}[t][t]{}%
\psfrag{x27}[t][t]{}%
\psfrag{x28}[t][t]{}%
\psfrag{x29}[t][t]{0.8}%
%
\psfrag{v01}[r][r]{0}%
\psfrag{v02}[r][r]{0.1}%
\psfrag{v03}[r][r]{0.2}%
\psfrag{v04}[r][r]{0.3}%
\psfrag{v05}[r][r]{0.4}%
\psfrag{v06}[r][r]{0.5}%
\psfrag{v07}[r][r]{0.6}%
\psfrag{v08}[r][r]{0.7}%
\psfrag{v09}[r][r]{0.8}%
\psfrag{v10}[r][r]{0.9}%
\psfrag{v11}[r][r]{1}%
\psfrag{v12}[r][r]{0}%
\psfrag{v13}[r][r]{}%
\psfrag{v14}[r][r]{}%
\psfrag{v15}[r][r]{}%
\psfrag{v16}[r][r]{}%
\psfrag{v17}[r][r]{}%
\psfrag{v18}[r][r]{}%
\psfrag{v19}[r][r]{}%
\psfrag{v20}[r][r]{0.8}%
\psfrag{v21}[r][r]{0}%
\psfrag{v22}[r][r]{}%
\psfrag{v23}[r][r]{}%
\psfrag{v24}[r][r]{}%
\psfrag{v25}[r][r]{}%
\psfrag{v26}[r][r]{}%
\psfrag{v27}[r][r]{}%
\psfrag{v28}[r][r]{}%
\psfrag{v29}[r][r]{0.8}%
%
\includegraphics[height=2.5in,width=6.5in,keepaspectratio]{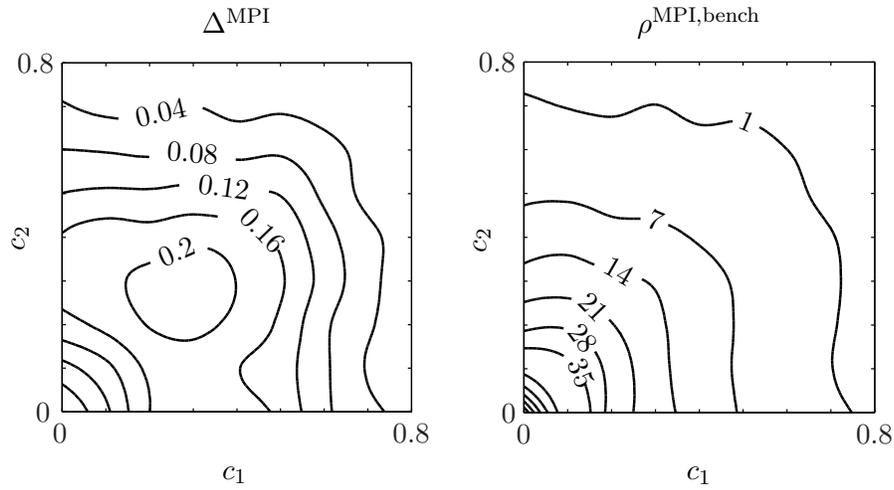}
\end{psfrags}%
\caption{Exp.\ 4: Average Relative (\%) Performance of MPI Policy vs.\ $(c_1, c_2)$, for
  $\beta = 0.9$.}
\label{fig:exp3plot}
\end{figure}

The fourth experiment assessed the effect of 
asymmetric constant startup costs  $(c_1, c_2) \in [0, 1]^2$ in two-project instances with no shutdown
costs and $\beta = 0.9$.
Figure  \ref{fig:exp3plot} displays the resultant contour plots for
$\Delta^{\textup{MPI}}$ and $\rho^{\textup{MPI}, \textup{bench}}$.
We obtain that $\Delta^{\textup{MPI}}$
reaches a maximum
value of about $0.23\%$, 
 vanishing as both startup costs approach zero 
and as either grows large enough. 
As for 
$\rho^{\textup{MPI}, \textup{bench}}$, it reaches a maximum value
of about $39.4\%$, and vanishes as either startup cost grows large enough and 
as both costs approach zero.

\begin{figure}[ht!]
\centering
\begin{psfrags}%
\psfragscanon%
\psfrag{s03}[t][t]{\color[rgb]{0,0,0}\setlength{\tabcolsep}{0pt}\begin{tabular}{c}$\beta$\end{tabular}}%
\psfrag{s04}[t][b]{\color[rgb]{0,0,0}\setlength{\tabcolsep}{0pt}\begin{tabular}{c}$\Delta^{\textup{MPI}}$\end{tabular}}%
\psfrag{s06}[t][t]{\color[rgb]{0,0,0}\setlength{\tabcolsep}{0pt}\begin{tabular}{c}$\beta$\end{tabular}}%
\psfrag{s07}[t][b]{\color[rgb]{0,0,0}\setlength{\tabcolsep}{0pt}\begin{tabular}{c}$\rho^{\textup{MPI}, \textup{bench}}$\end{tabular}}%
%
\psfrag{x01}[t][t]{0}%
\psfrag{x02}[t][t]{0.1}%
\psfrag{x03}[t][t]{0.2}%
\psfrag{x04}[t][t]{0.3}%
\psfrag{x05}[t][t]{0.4}%
\psfrag{x06}[t][t]{0.5}%
\psfrag{x07}[t][t]{0.6}%
\psfrag{x08}[t][t]{0.7}%
\psfrag{x09}[t][t]{0.8}%
\psfrag{x10}[t][t]{0.9}%
\psfrag{x11}[t][t]{1}%
\psfrag{x12}[t][r]{0.2}%
\psfrag{x13}[t][t]{}%
\psfrag{x14}[t][t]{}%
\psfrag{x15}[t][t]{}%
\psfrag{x16}[t][t]{}%
\psfrag{x17}[t][t]{}%
\psfrag{x18}[t][t]{}%
\psfrag{x19}[t][t]{0.9}%
\psfrag{x20}[t][r]{0.2}%
\psfrag{x21}[t][t]{}%
\psfrag{x22}[t][t]{}%
\psfrag{x23}[t][t]{}%
\psfrag{x24}[t][t]{}%
\psfrag{x25}[t][t]{}%
\psfrag{x26}[t][t]{}%
\psfrag{x27}[t][t]{0.9}%
%
\psfrag{v01}[r][r]{0}%
\psfrag{v02}[r][r]{0.1}%
\psfrag{v03}[r][r]{0.2}%
\psfrag{v04}[r][r]{0.3}%
\psfrag{v05}[r][r]{0.4}%
\psfrag{v06}[r][r]{0.5}%
\psfrag{v07}[r][r]{0.6}%
\psfrag{v08}[r][r]{0.7}%
\psfrag{v09}[r][r]{0.8}%
\psfrag{v10}[r][r]{0.9}%
\psfrag{v11}[r][r]{1}%
\psfrag{v12}[r][r]{0\%}%
\psfrag{v13}[r][r]{}%
\psfrag{v14}[r][r]{}%
\psfrag{v15}[r][r]{}%
\psfrag{v16}[r][r]{}%
\psfrag{v17}[r][r]{}%
\psfrag{v18}[r][r]{}%
\psfrag{v19}[r][r]{}%
\psfrag{v20}[r][r]{4\%}%
\psfrag{v21}[r][r]{0\%}%
\psfrag{v22}[r][r]{}%
\psfrag{v23}[r][r]{}%
\psfrag{v24}[r][r]{}%
\psfrag{v25}[r][r]{}%
\psfrag{v26}[r][r]{}%
\psfrag{v27}[r][r]{}%
\psfrag{v28}[r][r]{0.14\%}%
%
\includegraphics[height=2in,width=6.5in]{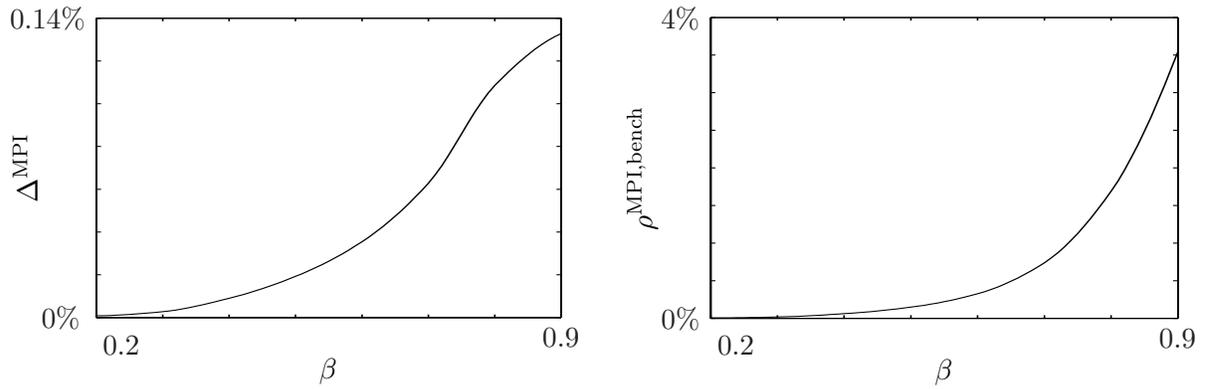}
\end{psfrags}%
\caption{Exp.\ 5: Average Performance of MPI policy for State-Dependent Startup Costs.}
\label{fig:exp4aplot}
\end{figure}

The fifth experiment evaluated the effect of 
state-dependent startup costs in two-project instances with no shutdown
costs as $\beta$ varies.
We  had MATLAB generate
 independent Uniform[0, 1] state-dependent startup costs for each
 instance. 
Figure \ref{fig:exp4aplot}'s  left pane plots $\Delta^{\textup{MPI}}$ vs.\ $\beta$, showing
the former to be increasing in the latter while
remaining well below $0.14\%$. 
The right pane 
plots 
$\rho^{\textup{MPI}, \textup{bench}}$,  which also increases
with
$\beta$, and remains below $4\%$.

\begin{figure}[ht!]
\centering
\begin{psfrags}%
\psfragscanon%
%
\psfrag{s01}[t][t]{\color[rgb]{0,0,0}\setlength{\tabcolsep}{0pt}\begin{tabular}{c}$c$\end{tabular}}%
\psfrag{s02}[t][b]{\color[rgb]{0,0,0}\setlength{\tabcolsep}{0pt}\begin{tabular}{c}$\Delta^{\textup{MPI}}$\end{tabular}}%
\psfrag{s03}[b][b]{\color[rgb]{0,0,0}\setlength{\tabcolsep}{0pt}\begin{tabular}{c}Dependence on $c$ for Multiple $\beta$\end{tabular}}%
\psfrag{s05}[t][t]{\color[rgb]{0,0,0}\setlength{\tabcolsep}{0pt}\begin{tabular}{c}$\beta$\end{tabular}}%
\psfrag{s06}[t][b]{\color[rgb]{0,0,0}\setlength{\tabcolsep}{0pt}\begin{tabular}{c}$\Delta^{\textup{MPI}}$\end{tabular}}%
\psfrag{s07}[b][b]{\color[rgb]{0,0,0}\setlength{\tabcolsep}{0pt}\begin{tabular}{c}Dependence on $\beta$ for $c = 0.3$\end{tabular}}%
%
\psfrag{x01}[t][t]{0}%
\psfrag{x02}[t][t]{0.1}%
\psfrag{x03}[t][t]{0.2}%
\psfrag{x04}[t][t]{0.3}%
\psfrag{x05}[t][t]{0.4}%
\psfrag{x06}[t][t]{0.5}%
\psfrag{x07}[t][t]{0.6}%
\psfrag{x08}[t][t]{0.7}%
\psfrag{x09}[t][t]{0.8}%
\psfrag{x10}[t][t]{0.9}%
\psfrag{x11}[t][t]{1}%
\psfrag{x12}[t][r]{$0.2$}%
\psfrag{x13}[t][t]{}%
\psfrag{x14}[t][t]{}%
\psfrag{x15}[t][t]{}%
\psfrag{x16}[t][t]{}%
\psfrag{x17}[t][t]{}%
\psfrag{x18}[t][t]{}%
\psfrag{x19}[t][t]{$0.9$}%
\psfrag{x20}[t][t]{$0$}%
\psfrag{x21}[t][t]{}%
\psfrag{x22}[t][t]{}%
\psfrag{x23}[t][t]{}%
\psfrag{x24}[t][t]{}%
\psfrag{x25}[t][t]{}%
\psfrag{x26}[t][t]{}%
\psfrag{x27}[t][t]{}%
\psfrag{x28}[t][t]{}%
\psfrag{x29}[t][t]{}%
\psfrag{x30}[t][t]{$1$}%
%
\psfrag{v01}[r][r]{0}%
\psfrag{v02}[r][r]{0.1}%
\psfrag{v03}[r][r]{0.2}%
\psfrag{v04}[r][r]{0.3}%
\psfrag{v05}[r][r]{0.4}%
\psfrag{v06}[r][r]{0.5}%
\psfrag{v07}[r][r]{0.6}%
\psfrag{v08}[r][r]{0.7}%
\psfrag{v09}[r][r]{0.8}%
\psfrag{v10}[r][r]{0.9}%
\psfrag{v11}[r][r]{1}%
\psfrag{v12}[r][r]{$0\%$}%
\psfrag{v13}[r][r]{}%
\psfrag{v14}[r][r]{}%
\psfrag{v15}[r][r]{}%
\psfrag{v16}[r][r]{}%
\psfrag{v17}[r][r]{}%
\psfrag{v18}[r][r]{$0.3\%$}%
\psfrag{v19}[r][r]{$0\%$}%
\psfrag{v20}[r][r]{}%
\psfrag{v21}[r][r]{}%
\psfrag{v22}[r][r]{}%
\psfrag{v23}[r][r]{}%
\psfrag{v24}[r][r]{}%
\psfrag{v25}[r][r]{$0.3\%$}%
%
\includegraphics[height=2in,width=6.5in]{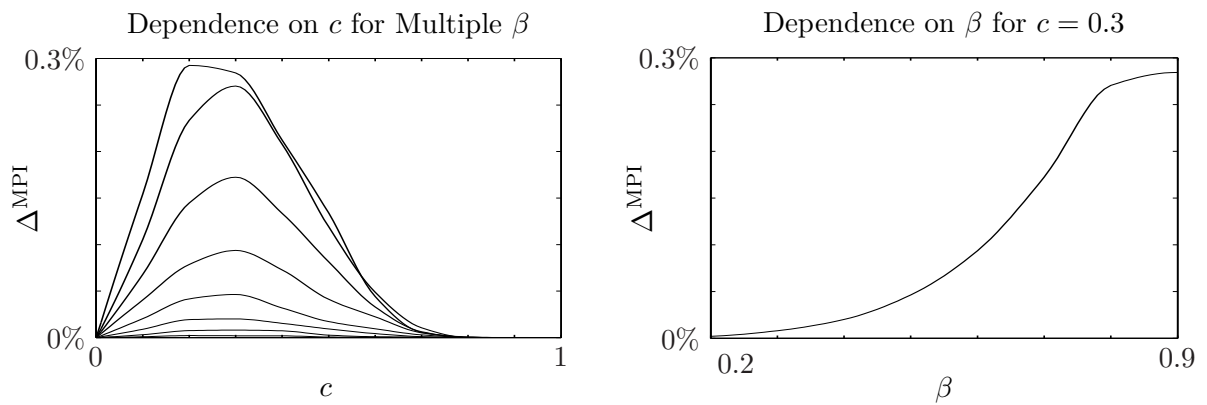}
\end{psfrags}%
\caption{Exp.\ 6: Counterpart of Figure \ref{fig:exp1aplot} for Three-Project Instances.}
\label{fig:plotThreeArmExp}
\end{figure}

\begin{figure}[ht!]
\centering
\begin{psfrags}%
\psfragscanon%
%
\psfrag{s01}[t][t]{\color[rgb]{0,0,0}\setlength{\tabcolsep}{0pt}\begin{tabular}{c}$c$\end{tabular}}%
\psfrag{s02}[b][b]{\color[rgb]{0,0,0}\setlength{\tabcolsep}{0pt}\begin{tabular}{c}$\rho^{\textup{MPI}, \textup{bench}}$\end{tabular}}%
\psfrag{s03}[b][b]{\color[rgb]{0,0,0}\setlength{\tabcolsep}{0pt}\begin{tabular}{c}Dependence on $c$ for Multiple $\beta$\end{tabular}}%
\psfrag{s05}[t][t]{\color[rgb]{0,0,0}\setlength{\tabcolsep}{0pt}\begin{tabular}{c}$c$\end{tabular}}%
\psfrag{s06}[b][b]{\color[rgb]{0,0,0}\setlength{\tabcolsep}{0pt}\begin{tabular}{c}$\rho^{\textup{MPI}, \textup{bench}}$\end{tabular}}%
\psfrag{s07}[b][b]{\color[rgb]{0,0,0}\setlength{\tabcolsep}{0pt}\begin{tabular}{c}Dependence on $\beta$ for $c = 0.1$\end{tabular}}%
%
\psfrag{x01}[t][t]{0}%
\psfrag{x02}[t][t]{0.1}%
\psfrag{x03}[t][t]{0.2}%
\psfrag{x04}[t][t]{0.3}%
\psfrag{x05}[t][t]{0.4}%
\psfrag{x06}[t][t]{0.5}%
\psfrag{x07}[t][t]{0.6}%
\psfrag{x08}[t][t]{0.7}%
\psfrag{x09}[t][t]{0.8}%
\psfrag{x10}[t][t]{0.9}%
\psfrag{x11}[t][t]{1}%
\psfrag{x12}[t][r]{$0.2$}%
\psfrag{x13}[t][t]{}%
\psfrag{x14}[t][t]{}%
\psfrag{x15}[t][t]{}%
\psfrag{x16}[t][t]{}%
\psfrag{x17}[t][t]{}%
\psfrag{x18}[t][t]{}%
\psfrag{x19}[t][t]{$0.9$}%
\psfrag{x20}[t][t]{$0$}%
\psfrag{x21}[t][t]{}%
\psfrag{x22}[t][t]{}%
\psfrag{x23}[t][t]{}%
\psfrag{x24}[t][t]{}%
\psfrag{x25}[t][t]{}%
\psfrag{x26}[t][t]{}%
\psfrag{x27}[t][t]{}%
\psfrag{x28}[t][t]{}%
\psfrag{x29}[t][t]{}%
\psfrag{x30}[t][t]{$1$}%
%
\psfrag{v01}[r][r]{0}%
\psfrag{v02}[r][r]{0.1}%
\psfrag{v03}[r][r]{0.2}%
\psfrag{v04}[r][r]{0.3}%
\psfrag{v05}[r][r]{0.4}%
\psfrag{v06}[r][r]{0.5}%
\psfrag{v07}[r][r]{0.6}%
\psfrag{v08}[r][r]{0.7}%
\psfrag{v09}[r][r]{0.8}%
\psfrag{v10}[r][r]{0.9}%
\psfrag{v11}[r][r]{1}%
\psfrag{v12}[r][r]{$0\%$}%
\psfrag{v13}[r][r]{}%
\psfrag{v14}[r][r]{}%
\psfrag{v15}[r][r]{}%
\psfrag{v16}[r][r]{}%
\psfrag{v17}[r][r]{}%
\psfrag{v18}[r][r]{}%
\psfrag{v19}[r][r]{}%
\psfrag{v20}[r][r]{$40\%$}%
\psfrag{v21}[r][r]{$0\%$}%
\psfrag{v22}[r][r]{}%
\psfrag{v23}[r][r]{}%
\psfrag{v24}[r][r]{}%
\psfrag{v25}[r][r]{}%
\psfrag{v26}[r][r]{}%
\psfrag{v27}[r][r]{}%
\psfrag{v28}[r][r]{}%
\psfrag{v29}[r][r]{$40\%$}%
\includegraphics[height=2in,width=6.5in]{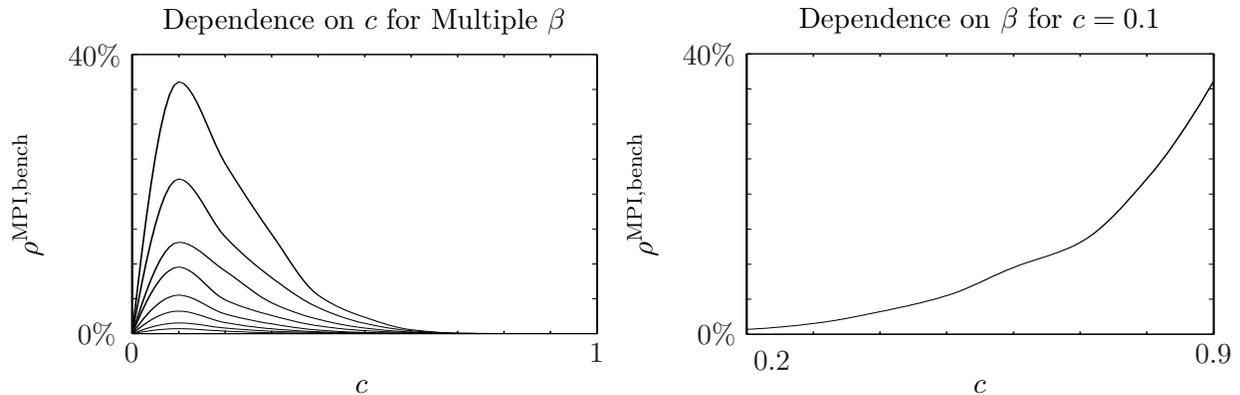}
\end{psfrags}%
\caption{Exp.\ 6: Counterpart of Figure \ref{fig:exp1bplot} for Three-Project Instances.}
\label{fig:plotThreeArmExpb}
\end{figure}

The sixth and last experiment is a three-project counterpart of the
second experiment,  based on a random sample of
100 instances of three 8-state projects each. 
The results are shown in Figures \ref{fig:plotThreeArmExp} and
\ref{fig:plotThreeArmExpb}, which are the counterparts of experiment 2's 
Figures \ref{fig:exp1aplot} and \ref{fig:exp1bplot}.
Comparison of Figures  \ref{fig:exp1aplot} and \ref{fig:plotThreeArmExp} 
reveals a slight performance degradation of the MPI policy in the
latter, though $\Delta^{\textup{MPI}}$ still remains quite small, 
below $0.3\%$. 
Comparison of Figures \ref{fig:exp1bplot} and \ref{fig:plotThreeArmExpb}
reveals similar values for  $\rho^{\textup{MPI},
  \textup{bench}}$.

\section{Conclusions}
\label{s:concl}
The paper has demonstrated the tractability and usefulness of the index
policy for bandits with switching costs based on the switching index introduced by \citet{asatene96}. 
The approach has been based on deploying the 
restless bandit indexation theory introduced by \citet{whit88} and
developed by the author  (cf.\ \citet{nmtop07}). 
\citet{nmswp07} announces results on the 
model extension that also
incorporates switching delays. 
The present analyses extend only in part part to such a case, as the
restless reformulation then yields semi-Markov projects that need
not be PCL-indexable. 
Other extensions of interest concern
restless bandits and particular models, as in
\citet{reiwe98}. 

\section*{Acknowledgments}
The author thanks the Associate Editor and two reviewers for
suggestions that helped improve the paper.
This work, which he started as a faculty member
at Universidad 
Pompeu Fabra (Barcelona), was
supported by the Spanish Ministry of Education \& Science
under projects BEC2000-1027 and MTM2004-02334, a Ram\'on y Cajal Investigator Award and an
I3 grant, by the European Union's Networks of Excellence 
EuroNGI and EuroFGI, and by the Autonomous Community of Madrid under grants
UC3M-MTM-05-075 and CCG06-UC3M/ESP-0767.
He presented part of this work at the Schloss Dagstuhl Seminar on
Algorithms for Optimization under Incomplete Information (Wadern,
Germany, 2005).


\end{document}